\documentclass[11pt,reqno]{amsart}
\usepackage{amsmath,amsthm,amsfonts,amssymb,mathrsfs,bm}
\usepackage[usenames]{color}
\usepackage{hyperref}
\usepackage{amssymb}
\usepackage{cases}
\usepackage{graphicx}
\usepackage{chngcntr}
\counterwithin{figure}{section}
\parskip        \bigskipamount

\newtheorem{theorem}{Theorem}[section]
\newtheorem{prediction}[theorem]{Prediction}
\newtheorem{lemma}[theorem]{Lemma}
\newtheorem{corollary}[theorem]{Corollary}
\newtheorem{proposition}[theorem]{Proposition}

\newtheorem{assumption}[theorem]{Assumption}

\newtheorem{remark}[theorem]{Remark}

\def\R{\mathbb{R}}

\def\E{\mathbb{E}}
\def\P{\mathbb{P}}

\def\cJ{\mathcal{J}}
\def\cW{\mathcal{W}}
\def\log{\mathrm{log}}

\def\min{\mathrm{min}}
\def\M{\mathcal{M}}
\def\S{\Sigma_{\mathrm{a}}}

\def\I{\mathcal{I}}
\def\T{\mathfrak{T}}

\def\supp{\mbox{\rm supp}}
\def\bx{{\bf x}}
\def\by{{\bf y}}

\def\bbx{\bar {\bf x}}
\def\bX{{\bf X}}
\def\bbX{\bar {\bf X}}

\def\bh{{\bf h}}
\def\bg{{\bf g}}
\def\bG{{\bf G}}
\def\bL{{\bf L}}
\def\bbh{\bar {\bf h}}
\def\bl{{\boldsymbol \lambda}}

\def\S{{\mathbb S}}
\def\MM{{M_1([\lambda_-^*,\lambda_+^*])}}

\makeatletter
\renewcommand*{\@cite@ofmt}{\hbox}

\newcommand{\abbr}[1]{{\sc\lowercase{#1}}}
\numberwithin{equation}{section}

\makeatother

\begin{document}

\author{Amir Dembo$*$
\and
Ofer Zeitouni$\dagger$
}
\thanks{$*$ Stanford University. Research partially supported by NSF grant DMS-1106627.}
\thanks{$\dagger$ Weizmann Institute of Science and Courant Institute.
Research partially supported by
a grant from the Israel Science Foundation.}
\date{\today}

\date{\today}
\title{Matrix optimization under random external fields}
\maketitle
 \begin{abstract}
We consider the quadratic optimization problem
$$F_n^{W,\bh}:=
\sup_{\bx \in S^{n-1}}
\Big( \frac12 \bx^T W \bx + \bh^T \bx \Big)\,,
$$
with $W$ a (random) matrix and $\bh$ a random external field. We study the
probabilities of large deviation of $F_n^{W,\bh}$ for $\bh$ a 
centered Gaussian vector with i.i.d. entries, both conditioned on 
$W$ (a general Wigner matrix), and unconditioned when $W$ 
is a \abbr{GOE} matrix. Our results validate (in a certain region) 
and correct (in another region), the prediction obtained 
by the mathematically non-rigorous replica method in
\textit{ Y. V. Fyodorov, P. Le Doussal,
J. Stat. phys. {\bf 154} (2014)}.
 \end{abstract}

\section{Introduction}
\label{sec-intro}

We consider in this paper the following quadratic optimization problem:
given an $n$-by-$n$ symmetric matrix $W$ and a vector $\bh\in \R^n$,
define, for 
$\bx\in S^{n-1}=\{\bx: \sum_{i=1}^n x_i^2=1\}$ in  the $n$-sphere, the 
quantity
\begin{equation}\label{eq:Enw-x-def}
E_n^{W,\bh}(\bx) :=
\frac12 \bx^T W \bx + \bh^T \bx\,,
\end{equation}
and consider the optimization problem
\begin{equation}\label{eq:Fnw-h-def}
F_n^{W,\bh}:=
\sup_{\bx \in S^{n-1}} \{ E_n^{W,\bh}(\bx) \} \,.
\end{equation}
As discussed
in \cite{FLD} to which we refer the reader for motivation and background,
the quantity $E_n^{W,\bh}(\bx)$ has a natural interpretation as minus the
energy 
associated with a configuration $\bx$ of $n$ spin variables $x_i$
in the presence of quadratic interaction $W$ and an external field
$\bh$. In contrast to the situation when $\bh=0$, the function
$F_n^{W,\bh}$ depends on the whole spectrum of $W$ and not just on its top
eigenvalue.

It is natural to consider both $W$ and $\bh$ as random objects.
Fixing $\Gamma>0$ constant, in \cite{FLD}, 
the authors (among other things) use a mathematically non-rigorous
replica method to study
the large deviations of the random variable
$F_n^{W,\bh}$ under the law $\P^{G,n}_\Gamma$ where
$\bh$ is a vector consisting of i.i.d. centered Gaussian variables of variance 
$\Gamma/n$ and $W$ is a matrix sampled from the 
Gaussian Orthogonal Ensemble (\abbr{GOE}). That is, $W$ 
is a symmetric matrix whose entries $\{W_{ij}\}$
on and above the
diagonal are independent centered Gaussian variables of variance 
$n^{-1}(1+\delta_{ij})$. In this setting, \cite{FLD} provides 
an argument for what
we refer to below as an \textit{annealed} Large Deviation Principle
(\abbr{LDP}), in the following form (see \cite[formula (43)]{FLD}).
\begin{prediction}
\label{theo-FLD}
Set $m_c:=\sqrt{1+\Gamma/(1+\Gamma)}$. Then, for any $m\in (m_c,\infty)$,
\begin{equation}
\label{eq-FLD1}
\lim_{\delta\to 0} \lim_{n\to\infty} 
\,\frac1n \log \, \P^{G,n}_{\Gamma}\,(\,|F_n^{W,\bh}-m|<\delta)=-I^{G,A}_{FLD}(m;\Gamma)\,,
\end{equation}
where
\begin{align}
I^{G,A}_{FLD}(m;\Gamma)=
\frac{m}{1+2\Gamma}\Big(-m\Gamma&+(1+\Gamma)\sqrt{m^2-m_c^2}\Big)
\notag\\
&-\log \left(\frac{\sqrt{1+\Gamma}
}{1+2\Gamma}\left(m+\sqrt{m^2-m_c^2}\right)\right)\,.
\label{eq-FLD2}
\end{align}
\end{prediction}
Note that no information is provided in Prediction \ref{theo-FLD}
on what happens when $m\leq m_c$.

As mentioned above, the derivation in \cite{FLD} uses a non-rigorous 
replica trick and breaks down at $m_c$. Our interest in the problem 
was initiated Y. Fyodorov, who asked whether Prediction \ref{theo-FLD} 
can be derived rigorously, 
and whether Prediction \ref{theo-FLD} can be extended
to the regime $m\leq m_c$.
This paper is devoted to answering these and related questions.

We find it advantageous and interesting
to discuss first a \textit{quenched}
large deviations theorem, 
namely a large deviations statement when the 
sequence of matrices $W=W_n$ is
given. In this setup, the assumption that $W$ is a \abbr{GOE} matrix
(or, more generally, a Wigner matrix)
plays no role.  Under appropriate conditions summarized in 
Assumption \ref{ass-QLDP}, we 
derive in Theorem 
\ref{theo-QLDP}
a conditional (in $\{W_n\}$) LDP 
for $F_n^{W,\bh}$ at speed $n$
when the vector $\bh$ is
either taken uniformly on $\sqrt{\Gamma}S^{n-1}$ (with associated 
explicit Good Rate Function (\abbr{GRF}) $I_q^H$), 
or
when the entries of $\bh$ are i.i.d. centered 
Gaussians with variance $\Gamma/n$ (with associated
\abbr{GRF} $I_q^G$).
(See for example \cite[Sec. 1.2]{DZ} for
definitions of \abbr{LDP} and \abbr{GRF}).
Theorem \ref{theo-QLDP} then yields in a straightforward manner Corollary  
\ref{cor-QLDP-W}, which deals with general Wigner matrices (including, 
but not limited to,  the
\abbr{GOE}).

We then turn our attention to
the case where $W$ 
is sampled from
the \abbr{GOE}. We derive the 
corresponding annealed (i.e. averaged on $W$) 
\abbr{LDP} at speed $n$ on the whole real line, see
Corollary \ref{cor-ALDP}. The proof builds on 
our quenched \abbr{LDP}, together with the \abbr{LDP} for the
top eigenvalue of Wigner matrices derived previously in \cite{BDG}. 

In the last subsection of the introduction, we 
simplify the general form of the quenched and annealed rate functions
for the \abbr{GOE}. 
In particular, 
we 
show that 
Prediction 
\ref{theo-FLD} is only true for $m>m_L:=1+\frac{\Gamma}{2(1+\Gamma)}$. Since 
$m_L>m_c$, this means that Prediction \ref{theo-FLD} does not hold in 
part of its domain (see Fig. \ref{fig-2} for a numerical example).

The annealed \abbr{LDP} for $F_n^{W,\bh}$
and our proof of it, are applicable 
more generally to any 
ensemble of random matrices having negligible
fluctuations of their empirical spectral measures 
at our large deviations speed and scale, and  
for which the \abbr{LDP} at speed $n$ of the 
maximal (or minimal), eigenvalue 
of $W$ is available. In contrast with the 
\emph{universality} of the rate functions $I^H$ 
and $I^G$ for the \emph{quenched} large deviations 
of $F_n^{W,\bh}$, the annealed rate functions
$I^{H,A}$ and $I^{G,A}$ are specific to the \abbr{GOE}
(as they depend on the exact form of the \abbr{LDP}
for its maximal eigenvalue). 
 
In the rest of the introduction we present the relevant notation 
and state our assumptions and main results. 

\label{intro}
\subsection{LDP for quadratic optimization problems}
Throughout we write $\bx=(x_1,\ldots,x_n)$ 
for a vector in $\R^n$ and
$\bx^2=(x_1^2,\ldots,x_n^2)$.
The scalar product in $\R^n$ is denoted $\langle\cdot,\cdot\rangle$, with $\|\cdot\|$ for the Euclidean norm in $\R^n$. We further use 
$M_+(J)$ for the space of all finite, Borel measures on $J \subseteq \R$,
with $M_1(J)$ denoting the sub-space of all probability measures on $J$, both 
equipped with the topology of weak convergence.

Let $\R^n_{\geq}:=\{(\lambda_1,\lambda_2,\ldots,\lambda_n)\in \R^n:
\lambda_1\geq \lambda_2\geq\ldots\geq \lambda_n\}$ denote the collection
of ordered $n$-tuple real numbers, with $S^{n-1}$ denoting
the usual Euclidean sphere of radius $1$.
For fixed $\bl\in \R^n_{\geq}$ and constant $\Gamma > 0$, 
we are interested in large deviations for 
the (random) optimization problem
\begin{equation}
        \label{eq-opt}
        F^*_{n,\bh,\bl}=\sup_{ \bx \in \, S^{n-1}} 
        \left(
\frac12 \langle \bl,\bx^2\rangle+  
\langle \bh,\bx\rangle\right)\,,
\end{equation}
with respect to $\bh=(h_1,\ldots,h_n)\in \sqrt{\Gamma} S^{n-1}$ 
a random vector drawn uniformly from the Haar measure on 
$\sqrt{\Gamma} S^{n-1}$,
or alternatively when having $\bh=\bg$ a
centered multivariate normal random vector  
of covariance matrix $\frac{\Gamma}{n} {\bf I}_n$.
Throughout we assume the following about the 
parameters of the optimization problem \eqref{eq-opt}.
\begin{assumption}
        \label{ass-QLDP}
        For $n \to \infty$ we have that:
\newline        
(A1). 
$L_n^{\bl}=\frac1n \sum_{j=1}^n \delta_{\lambda_j}$ 
converge weakly in $M_1(\R)$ to some $q(\cdot)$ 
of compact support.
\newline
(A2). $\lambda_1(n)\to \lambda_+^*<\infty$
(necessarily, $\lambda_+^* \ge 
q_+ := 
\max\{x\in \supp(q)\}$).
\newline
(A3).
$\lambda_n(n)\to \lambda_-^*>-\infty$
(necessarily, $\lambda_-^* \le q_-:=
\min\{x\in \supp(q)\}$).                                      
\end{assumption}

Our first result is then the following \abbr{LDP}.
\begin{theorem}
        \label{theo-QLDP}
        Let Assumption \ref{ass-QLDP} hold
        and fix $\Gamma > 0$ non-random.               
\newline 
(a). For 
$\bh$ Haar distributed on $\sqrt{\Gamma} S^{n-1}$,
        the sequence $\{F^*_{n,\bh,\bl}\}$ satisfies 
        the \abbr{LDP} in $\R$,
        with speed $n$ and \abbr{GRF}
        \begin{equation}
                \label{eq-rateQ}
                I^H_q(m;\lambda^*_\pm,\Gamma)=
                        \inf\left\{ \frac12 H(q|\nu): 
                        \nu \in M_1([\lambda_-^*,\lambda_+^*]), 
                        m=F(\lambda_+^*,\nu;\Gamma)\right\},
                        \end{equation}
where for given $\xi \in \R$,
\begin{equation}\label{eq:Mnu-def}
F(\xi,\nu;\Gamma) = 
\frac12 \inf_{\theta> \xi}\left[\theta+ \Gamma \int\frac{\nu(dx)}{\theta-x} \right]\,,
\end{equation}
for any $\nu \in M_+((-\infty,\xi])$, and
$$
H(\mu|\nu) =\left\{\begin{array}{ll} \int d\mu \log (\frac{d\mu}{d\nu})
+ \nu(\R) - 1, & \mu \ll \nu,\\
\infty, & \mbox{\rm otherwise}.
\end{array}\right. 
$$
(b). For $\bg$ centered 
multivariate normal of covariance $\frac{\Gamma}{n} {\bf I}_n$, 
the sequence  
$\{F^*_{n,\bg,\bl}\}$ satisfies the \abbr{LDP} with speed $n$ and 
the \abbr{GRF}
\begin{align}\label{eq:rf-gauss}
I^G_q(m;\lambda_\pm^*,\Gamma) &= 
\inf\left\{ \frac12 H(q|\nu): 
                        \nu \in M_+([\lambda_-^*,\lambda_+^*]), 
                        m=F(\lambda_+^*,\nu;\Gamma)\right\}\nonumber \\
&=
\inf_{y \ge 0} \big\{ I^H_q(m;\lambda_\pm^*,\Gamma y) + J_1 (y) \big\}
\,,
\end{align}
where 
\begin{equation}\label{eq:J-def}
J_1 (y) = \frac12 \Big(y -1 -\log \, y\, \Big) \,, \qquad \forall y \in \R_+ \,.  
\end{equation}
\end{theorem}

\begin{remark}\label{rmk-mstar}
Clearly,
$I^G_q(\cdot;\lambda^*_\pm,\Gamma) \le 
I^H_q(\cdot;\lambda^*_\pm,\Gamma)$ and both
of these \abbr{GRF}s are zero if and only if 
$\nu=q$ and $y=1$. That is, when $m=\overline{m}$, where
\begin{equation}
\label{eq-mstar}
\overline{m} := F(\lambda_+^*,q;\Gamma) \,.
\end{equation}
Further, from \eqref{eq:Mnu-def} we see that 
$F(\lambda^*_+,\nu;\Gamma) \in [m^*_-,m^*_+]$
for $\nu \in M_1([\lambda^*_-,\lambda^*_+])$ and
\begin{align}\label{eq-m-}
m^{*}_{-} &:=
F(\lambda^{*}_{+},
\delta_{\lambda^{*}_{-}}
;\Gamma) =
\frac12 \big(        \theta^*_- + \frac{\Gamma}{\theta^*_--\lambda_-^*}\big)\,,
        \quad \text{where}\; \theta^*_-=\lambda_+^*\vee (\lambda_-^*+\sqrt{\Gamma})\,
\\
\label{eq-m+}
m_{+}^{*} &:= 
F(\lambda^{*}_{+},
\delta_{\lambda^{*}_{+}}
;\Gamma) =
         \frac12\lambda_+^*+ \sqrt{\Gamma}\,.
\end{align}        
Hence, $I^H_{q}(\cdot;\lambda^*_\pm,\Gamma)=\infty$ 
outside the compact interval $[m^*_-,m^*_+]$ (which is strictly 
above $\frac12 \lambda_+^*$), whereas 
$I^G_{q}(m;\lambda^*_\pm,\Gamma)=\infty$ for $m \le \frac12 \lambda_+^*$.
\end{remark}

We next detail a few regularity properties of the rate functions of 
Theorem \ref{theo-QLDP}.
\begin{proposition}\label{prop:rf-prop}
The \abbr{GRF} $I^H_{q}(\cdot;\lambda^*_\pm,\Gamma)$ 
is continuous on $(m_-^*,m_+^*)$, non-increasing
on $(m_-^*,\overline{m}]$ and convex strictly 
increasing on $[\overline{m},m_+^*)$, whereas
the \abbr{GRF} $I^G_q(\cdot;\lambda^*_\pm,\Gamma)$ 
is continuous on $(\frac{1}{2}\lambda^*_+,\infty)$,
non-increasing on 
$(\frac{1}{2}\lambda^*_+,\overline{m}]$ and 
convex strictly 
increasing on $[\overline{m},\infty)$.
\end{proposition}

\begin{remark}\label{rem:IG-lminus}
See also Proposition \ref{expl-rf} for more explicit expressions for the 
rate functions $I^H_q$ and $I^G_q$. In particular, it is shown there that  
$I^G_q(m;\lambda^*_{\pm},\Gamma)$ is \emph{independent} of $\lambda^*_-$,
as is $I^H_q(m;\lambda^*_{\pm},\Gamma)$ when $m > \overline{m}$, whereas
$I^H_q(m;\lambda^*_{\pm},\Gamma)$
is independent of $\lambda^*_+$ when $m < \overline{m}$. 
\end{remark}

\subsection{LDP for random quadratic forms - 
Wigner matrices versus the GOE}
The general \abbr{LDP} of Theorem \ref{theo-QLDP} yields \abbr{LDP}s
for quadratic optimization problems involving random matrices. 
Fixing $\lambda^*_\pm \in \R$ and 
$q \in M_1([\lambda^*_-,\lambda^*_+])$,
let $\cW_{\lambda^*_\pm,q}$ denote the collection 
of all sequences of (random or deterministic) 
symmetric $n$-dimensional $\R$-valued matrices, 
whose ordered eigenvalue vectors $\bl$ satisfy 
Assumption \ref{ass-QLDP} for $q$ and $\lambda^*_\pm$. 
The following \abbr{LDP} for 
$\{F_n^{W,\bh}\}$ is a direct consequence of 
Theorem \ref{theo-QLDP}. Here and in the sequel, for a sequence 
$\{W_n,\bh_n\}$ we write 
$F_n^{W,\bh}$ as shorthand for $F_n^{W_n,\bh_n}$.
\begin{corollary}[Quenched \abbr{LDP}]
\label{cor-QLDP-W}
Fix a deterministic constant $\Gamma > 0$   and
a sequence $\{W_n\} \in \cW_{\lambda^*_\pm,q}$.
For $\widetilde{\bh}$ independent of $W_n$, denote by
$\P^{W,H,n}_\Gamma$ the law of $F_n^{W,\widetilde{\bh}}$ conditioned 
on $W_n$ if $\widetilde{\bh}$ is Haar distributed on
$\sqrt{\Gamma}S^{n-1}$, and by $\P^{W,G,n}_\Gamma$ 
if $\widetilde{\bh}$ is a centered
multivariate normal $\widetilde{\bg}$, of covariance $\frac{\Gamma}{n}
{\bf I}_n$.
\newline
(a).  The sequence $\big\{\P^{W,H,n}_\Gamma \big\}_{n\geq 1}$  satisfies the 
\abbr{LDP} on $\R$ with  speed $n$ and the \abbr{GRF}
$I^H_{q}(m;\lambda^*_\pm,\Gamma)$ given
by \eqref{eq-rateQ} (or alternatively, 
\eqref{eq:const4}).
\newline                 
(b). 
The sequence 
$\big\{\P^{W,G,n}_{\Gamma}\big\}_{n \ge 1}$ satisfies the corresponding 
\abbr{LDP} with 
\abbr{GRF}
$I^G_q(m;\lambda^*_\pm,\Gamma)$
given by \eqref{eq:rf-gauss} (or alternatively, \eqref{eq:const3}). 
\end{corollary}

\begin{remark} Recall that a symmetric random 
matrix $W_n$ is a \textit{Wigner} matrix if it
has
centered independent entries
on and above the diagonal, with the entries above 
the diagonal being i.i.d. of variance $1$ and bounded fourth moment,
while the on-diagonal entries are i.i.d. with
uniformly bounded second moment.
Such matrices 
are a.s. in $\cW_{\pm 2,\sigma}$, 
with $\sigma$ the \textit{semi-circle} law having the
support $[-2,2]$ and density 
$f_\sigma(x)=(2\pi)^{-1} \sqrt{4-x^2} 1_{|x| \le 2}$, see
\cite[Theorem 2.1.21]{AGZ} and \cite{BY}.
Hence, all the conclusions of Corollary \ref{cor-QLDP-W} 
hold for such matrices.
\end{remark}

We turn to the \abbr{LDP} averaged over the choice of $W_n$ 
from the \abbr{GOE}.
Let $\P^{H,n}_\Gamma=
\E_{\mbox{\rm \abbr{GOE}}}\P^{W,H,n}_\Gamma$ and 
$\P^{G,n}_\Gamma= 
\E_{\mbox{\rm \abbr{GOE}}}\P^{W,G,n}_\Gamma$.
Note that under either $\P_\Gamma^{H,n}$ or $\P^{G,n}_\Gamma$, 
the matrix $W$ is sampled according to the \abbr{GOE} and is independent
of the random vector $\widetilde{\bh}$.
\begin{corollary}[Annealed \abbr{LDP}]
\label{cor-ALDP}
%
$~$
\newline
(a). 
The sequence $\big\{\P^{H,n}_\Gamma\big\}_{n\geq 1}$  satisfies the 
\abbr{LDP} with speed $n$ and \abbr{GRF}
\begin{equation}\label{eq:ann-rf}
I^{H,A}(m;\Gamma)=
\inf_{\psi^*_- \le -2, \psi^*_+ \ge 2}
\Big\{ I^H_{\sigma}(m;\psi^*_\pm,\Gamma) + 
I_e(\psi^*_+) + I_e(-\psi^*_-) \Big\},
\end{equation}
for $I^H(\cdot)$ of \eqref{eq-rateQ} (or alternatively \eqref{eq:const4}), and 
\begin{numcases}{\!\!\!\!\!\!\!\!\!\!\!\!
\label{eq:Ie-def}
I_e(\psi)=}
\int_2^\psi \sqrt{(u/2)^2-1} \, du, & $\psi \ge 2$, \\
\qquad \qquad \infty,& \textrm{otherwise}. \nonumber
\end{numcases}
\newline
(b). 
The sequence $\big\{\P^{G,n}_\Gamma \big\}_{n \ge 1}$ 
satisfies the \abbr{LDP} with speed $n$ and \abbr{GRF}
\begin{equation}\label{eq:ann-rf-gauss}
I^{G,A} (m;\Gamma) = \inf_{\psi^* \ge 2} \Big\{ 
I^G_{\sigma}(m;\pm \psi^*,\Gamma) 
+  I_{e} (\psi^*) \Big\},
\end{equation}
for $I^G(\cdot)$ of \eqref{eq:rf-gauss}
(or alternatively \eqref{eq:const3}).
\end{corollary}

\subsection{Explicit rate functions}
We shall derive explicit expressions for the various 
rate functions introduced in the article, starting with the
\abbr{GRF}-s of Theorem \ref{theo-QLDP}, for general $q \in M_1([\lambda^*_-,\lambda^*_+])$. To state the result,
we require the logarithmic potential and Stieltjes transform of $q(\cdot)$, denoted by
\begin{align}\label{eq:bL-def}
\bL(\xi) &= \int \log|\xi-x| q(dx) \,,\quad
\forall \xi \notin (\lambda^*_-,\lambda^*_+) \,,
 \\ 
\bG(\xi) = \bL'(\xi) &=\int (\xi-x)^{-1} q(dx)  \,,\quad
\forall \xi \notin (\lambda^*_-,\lambda^*_+) \,.
\label{eq:bG-def}
\end{align}

\begin{proposition}[Quenched rate functions]\label{expl-rf}
$~$
\newline
(a) In case $q_{\pm}=\lambda^*_{\pm}$, 
the \abbr{GRF} for part (a) of Theorem 
\ref{theo-QLDP} is
\begin{equation}\label{eq:const4}
I^H_q(m;\lambda^*_\pm,\Gamma) = \frac12 
\big[ \log |B| + \bL(\psi) - \bL(\theta) \big]\,.
\end{equation}
Here $t \ge 0$, $\theta \ge \lambda^*_+$,
$B$ and $\psi$ are such that
\begin{align}
\label{eq:const5}
B-1 &= \qquad \; (\theta-\psi)\bG(\psi) + Bt 
\\
B(2m-\theta) &= \qquad \qquad  \Gamma \bG(\psi) 
+ \frac{B \Gamma t}{\theta-\psi^*}   
\;\; \text{ with } t=0 \text{ whenever } \;
\psi \ne \psi^*\,,
\label{eq:const2}
\\
B(\psi-\theta) &\ge \Gamma \bG(\theta) - 
\Gamma \bG(\psi) - \frac{B \Gamma t}{\theta-\psi^*}
\; \text{ with equality whenever} \; 
\theta>\lambda^*_+\,,
\label{eq:const1}
\end{align}
with $m > \overline{m}$ requiring $B>0$ and 
$\theta > \psi \ge \psi^*=\lambda^*_+$,
while for $m < \overline{m}$ we consider
either $B>0$ and $\psi > \theta$, or 
$B<0$ and $\psi \le \psi^*=\lambda^*_-$.
\newline
(b) For any $q \in M_1([\lambda^*_-,\lambda^*_+])$, 
the \abbr{GRF} for part (b) of Theorem 
\ref{theo-QLDP} is
\begin{equation}\label{eq:const3}
I^G_q(m;\lambda^*_\pm,\Gamma) = 
\frac12 \big[ (\theta-\psi) \bG(\psi) + t + \bL(\psi) 
- \bL(\theta) \big]\,,
\end{equation}
with $B=1$ and $\theta,\psi \ge \psi^* = \lambda_+^*$,
$t \ge 0$, are determined by 
\eqref{eq:const2} and \eqref{eq:const1}.
\end{proposition}

Building on Proposition \ref{expl-rf}, due to the 
simple form of $\bG(\cdot)$ for the 
semi-circle, one can explicitly solve 
the variational problems in the definition of the rate functions 
\begin{equation}\label{eq:Q-ann-comp}
I^{G,A} \le I^{H,A} \wedge I^{G}_{\sigma}
\le I^{H,A} \vee I^{G}_{\sigma}
\le I^{H}_{\sigma} \,.
\end{equation}
To state the result, we introduce $\alpha,\beta \ge 1$
such that $\theta=\alpha+\alpha^{-1}$ and $\psi=\beta+\beta^{-1}$,
the functions
\begin{equation}\label{eq:T-q-def}
\I(\alpha,\beta)=
J_1\Big(\frac{\alpha}{\beta}\Big)
- \frac{1}{4}(\alpha^{-1}-\beta^{-1})^2\,,\quad
\T(\alpha)=
\Gamma^{-1} (\alpha+\frac{1}{\alpha}-2)[
2m-\alpha-\frac{1}{\alpha}-\Gamma] \,,
\end{equation}
and the constants
$1 < m_c < m_L < \overline{m} < m_U$, given by
\begin{equation}\label{eq:m-phases}
m_c=\sqrt{\frac{1+2\Gamma}{1+\Gamma}}\,,\quad
m_L=1+\frac{\Gamma}{2(1+\Gamma)}\,,\quad
\overline{m}=\sqrt{1+\Gamma}\,,\quad
m_U=1+\frac{\Gamma(1+2\Gamma)}{2(1+\Gamma)}\,.
\end{equation}
\begin{figure}[htb]
\vskip -3cm
\begin{center}
\includegraphics[width=120mm]{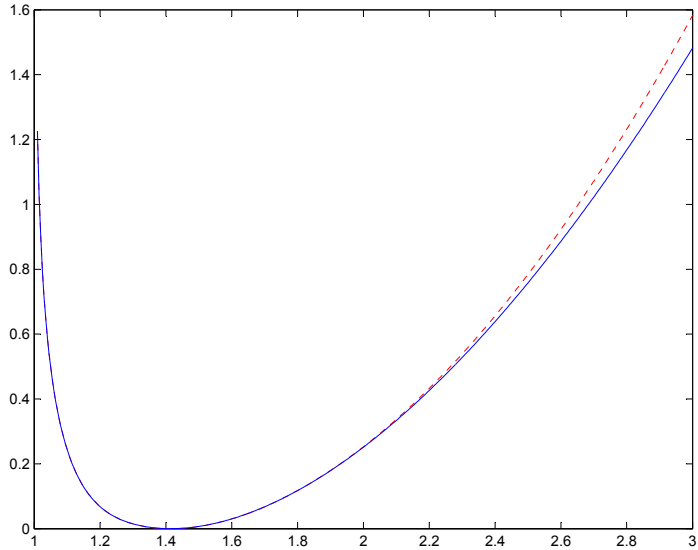}
\caption{The rate functions $I^{G,A}(m;1)<I^G_\sigma(m;\pm2,1)$
\footnotesize{
(here $m_L=1.25$, $m_U=1.75$)
}
}
	\label{fig-1}
\end{center}\end{figure}
\begin{figure}[htb]
\vskip -3cm
\begin{center}
\includegraphics[width=100mm]{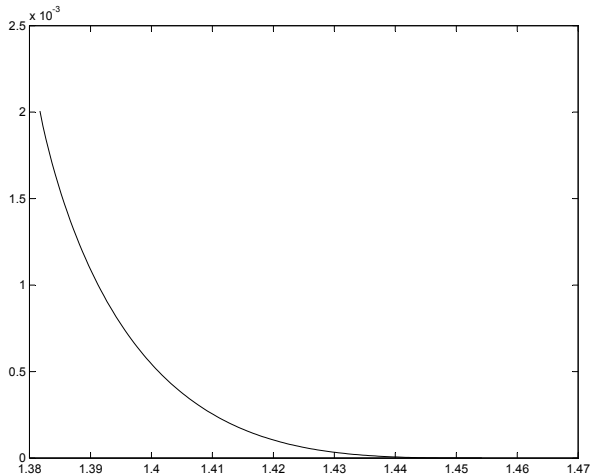}
\caption{
$I^{G,A}_{FLD}(m;10)-I^{G,A}(m;10)$ for $m\in (m_c,m_L)$;
\footnotesize{here $m_c=1.38[2]$, $m_L=1.45[4]$}
}
	\label{fig-2}
\end{center}\end{figure}
\begin{proposition}[Rate functions with semi-circle]\label{expl-rf-W}
$~$
\newline
(a). The quenched \abbr{GRF} $I_\sigma^G(\cdot;\pm2,\Gamma)$ is given by the formula
\begin{equation}\label{eq:rf-gauss-semi}
I^{G}_{\sigma}(m;\pm 2,\Gamma) = \I(\alpha_q,\beta_q)+\frac12 \T(\alpha_q)
{\bf 1}_{\{m \ge m_U\}}\,,
\end{equation}
where $\alpha_q \le \beta_q$ iff $m \le \overline{m}$ are
given by 
\begin{numcases}{\!\!\!\!\!\!\!\!\!\!\!\!
(\alpha_q,\beta_q)=}
\label{eq:theta-is-2} 
(1,\frac{\Gamma}{2(m-1)}),&
$m \in (1,m_L]$, \\
\label{eq:int-pt}
(m_c^{-2}[m+\sqrt{m^2-m_c^2}],(1+\Gamma)[m-\sqrt{m^2-m_c^2}]),&
$m \in (m_L,m_U), $ \\
\label{eq:psi-is-2}
( \frac12[(m+1)+\sqrt{(m+1)^2 -4-2\Gamma}],1),&  
$m \in [m_U,\infty). $
\end{numcases}
\newline
(b). The strictly convex annealed \abbr{GRF} $I^{G,A}(m;\Gamma)$
equals the quenched \abbr{GRF} 
$I^{G}_\sigma(m;\pm 2,\Gamma)$ from \eqref{eq:rf-gauss-semi}
for $m \in [1,m_U]$, whereas for $m > m_U$, 
\begin{align}
I^{G,A}(m;\Gamma)& =\I(\alpha_a,\beta_a^{-1})
\label{eq:ann-rf-gauss-exp}
\end{align}
for $(\alpha_a,\beta_a^{-1})$ as in \eqref{eq:int-pt}, i.e.
\begin{equation}
\label{eq-addendum}
(\alpha_a,\beta_a^{-1})=(m_c^{-2}[m+\sqrt{m^2-m_c^2}],(1+\Gamma)[m-
\sqrt{m^2-m_c^2}])\,. 
\end{equation}
\end{proposition}

See Figure \ref{fig-1} for a plot of the quenched and 
annealed rate functions $I^{G}_{\sigma}(m;\pm 2,1)$ and $I^{G,A}(m;1)$.
\begin{remark}
The annealed \abbr{GRF} $I^{G,A}(m;\Gamma)$ could also be written as 
$$
I^{G,A}(m;\Gamma)=\I(\alpha,\beta)
$$
where $(\alpha,\beta)=(\alpha_q,\beta_q)$ if $m\leq m_U$,
and $(\alpha,\beta)$ given by the r.h.s. of
\eqref{eq-addendum}, if $m>m_U$. Note that while
$m \mapsto I^{G,A}(m;\Gamma)$ is smooth except 
for the jump discontinuity of its third derivative 
at $m=m_L$, 
the function 
$m \mapsto I^G_\sigma(m;\pm 2, \Gamma)$ 
is also non-smooth at $m=m_U$.
\end{remark}

\begin{remark}
\label{rem-fyod}
It is worthwhile to comment on the relation between the rate function
$I^{G,A}(m;\Gamma)$ of Proposition \ref{expl-rf-W} and Prediction
\ref{theo-FLD} from \cite{FLD}: a tedious, but straight forward 
algebraic manipulation shows that 
$I^{G,A}(m;\Gamma)=I^{G,A}_{FLD}(m;\Gamma)$ for $m\geq m_L$. However, 
both a numerical evaluation, see Figure \ref{fig-2}, and analytic
evaluation of the limit $m\searrow m_c$
as well as comparison of
the first three derivatives at $m=m_L$,  
show that in general $I^{G,A}(m;\Gamma)\neq I^{G,A}_{FLD}(m;\Gamma)$ 
in the interval $m\in(m_c,m_L)$.
\end{remark}

\begin{figure}
 \begin{center}
\includegraphics[width=80mm]{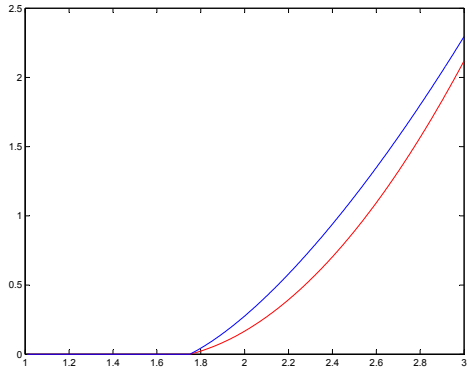}
\vspace{3mm}
\caption{The parameter $t(m)$ 
	\footnotesize{(quenched, dashed black; annealed, doted green)}, at $\Gamma=1$.}
	\label{fig-4}
\end{center}\end{figure}

We conclude the introduction with comments on the rate functions
$I^G_\sigma$ and $I^{G,A}$. The general form of the rate function can be
understood  by considering the following heuristics. 
Three main objects 
enter 
the (diagonalized) 
optimization problem \eqref{eq-opt}:
\begin{enumerate}
	\item The total mass $\sum_{i=1}^n h_i^2$, which we take to
roughly equal $\Gamma y$.
\item The measure of total mass $y>0$, which controls the 
      distribution of the $h_i^2$  
      as weights on the eigenvalues $\lambda_i$, denoted
	$\nu=\Gamma^{-1}\sum_{i=1}^n h_i^2 \delta_{\lambda_i}$.
\item The optimal profile of $x_i$-s for a given $\nu$, which turns out to
	be determined by a Lagrange multiplier 
	$\theta\geq \lambda_+^*$ (specifically, $x_i=h_i/(\theta-\lambda_i)$,
	as shown in Lemma \ref{lem-Jf}).
	\end{enumerate}
	The minimization of the
	probabilistic cost of producing such $\nu$ while 
	constraining the value $F^*_{n,\bh,\bl} \approx m$ yields for $n \to \infty$ 
	the optimal 
	\begin{equation}\label{eq:opt-nu}
	\nu^*(d\lambda)=\frac{\theta-\lambda}{\psi-\lambda}q(d\lambda) 
	+ t \delta_{\lambda^*_+} 
	\end{equation}
	in terms of another Lagrange multiplier, denoted $\psi\geq
	\lambda_+^*$ (see proof of 
	part (b) of Proposition \ref{expl-rf} for the derivation of $\nu^*$).	
    \begin{figure}    \begin{center}
\includegraphics[width=80mm]{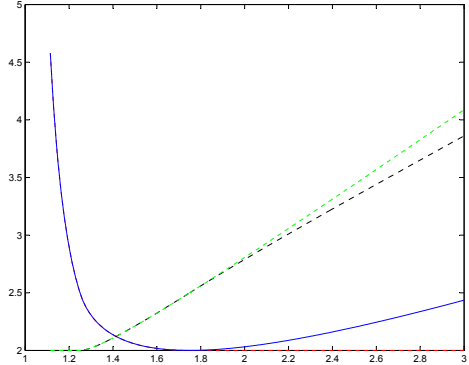}
\vspace{3mm}
\caption{The parameters $\psi(m)$ \footnotesize{(quenched, dashed red; 
	annealed, solid blue)} and
	$\theta(m)$ \footnotesize{(quenched, dashed black; annealed, doted green)}, 
	for $\Gamma=1$.}
	\label{fig-3}
\end{center}\end{figure}
	We note in passing that $t=t(m)$ represents the total 
	mass projected by $\bh$ on the eigenspace of $o(n)$ top 
    eigenvalues of $W$, if constrained to $F^{W,\bh}_n \approx m$,
    and is non-zero only when $\psi$ is at the edge of the vector $\bl$. 
	Now the
	three regimes of the quenched rate function in Proposition
    \ref{expl-rf-W}, where $\lambda^*_+=2$, correspond to the following cases:
	\begin{eqnarray*}
	m\in \big(1,m_L\big]     &\quad \Longleftrightarrow \quad&
	\psi> \theta=\lambda_+^*, \quad \quad \; t=0 \\
	m\in (m_L,m_U)   &\quad \Longleftrightarrow \quad&\theta>\lambda_+^*,\; \psi>\lambda_+^*, \; t=0\\
	m\in [m_U,\infty)&\quad \Longleftrightarrow \quad&\theta > \psi=\lambda_+^*, \qquad \; t>0\,.
	\end{eqnarray*}
	At the typical value $\overline{m}$ 
	which lies in 
	$(m_L,m_U)$, one switches from having $\psi>\theta$ (lower tail large deviations,
	with $y=y(m)<1$), to $\theta>\psi$ (upper tail large deviations, with $y=y(m)>1$).
	In the annealed case described in
	Proposition \ref{expl-rf-W}, the regime $m\in[m_U,\infty)$ is different
    because while saturating the constraint on $\psi$ at the value of 
    the top eigenvalue, the optimal solution is now able to shift the
    top eigenvalue from $\lambda^*_+$ to $\psi=\psi(m)>\lambda^*_+$. 
    See Figures \ref{fig-3}-\ref{fig-4} for a 
    plot of the parameters of \eqref{eq:opt-nu}
    in the quenched and annealed cases (at $\Gamma=1$, with Fig. \ref{fig-3}
    depicting $m \mapsto \psi(m), m\mapsto \theta(m)$, 
    and Fig. \ref{fig-4} for $m \mapsto t(m)$). 
\begin{remark} The parameter $\Gamma$ determines the relevant importance of the 
quadratic and linear parts of the optimization problem \eqref{eq-opt}. In one
extreme of $\Gamma \gg 1$, the typical value is $\overline{m} = \sqrt{\Gamma} (1+o(1))$
with \eqref{eq-opt} dominated by its linear (Gaussian) part. In this case,
except for the extreme tails (i.e. $m \le m_L \approx 3/2$ 
or $m \ge m_U \approx \Gamma$), both  
$I^{G,A}(m;\Gamma)$ and $I^G_\sigma(m;\pm 2, \Gamma)$ approximately match
the $\chi$-square rate $J_1(y)$ for $y=m^2/\Gamma$ 
(
with the contribution of eigenvalues buried in the correction terms). 
In contrast, for $\Gamma \downarrow 0$ the quadratic part dominates. Further,  
with $m_U(\Gamma) \downarrow 1$, its all about the top of the spectrum of $W_n$.
Here $I^{G,A}(m;\Gamma)$ approximately matches the \abbr{GOE} rate function 
$I_e(2m)$ (and since $\bh$ does not matter much, one uses $t(m)$ small 
and $\theta(m) \approx \psi(m) \approx 2m$ to get there), whereas 
$I^G_\sigma(m;\pm 2,\Gamma) \approx t/2 = (m-1)^2/(2\Gamma)$ is the 
cost of making up the
$m-1$ discrepency in the value of $F^{W,\bh}_n$ by having $\bh$ of that 
magnitude, aligned to the top eigenvector of $W_n$.
\end{remark}

\section{Proofs}

\subsection{Rate functions: regularity properties}

\proof[Proof of Proposition \ref{prop:rf-prop}]
Fixing $\lambda_{\pm}^*$ and $\Gamma$, let 
$\M=M_+([\lambda^*_-,\lambda^*_+])$, 
$\M_1:=\MM$,  and
$f(\theta,\nu):=\frac12\big(\theta+\Gamma\int (\theta-x)^{-1} \nu(dx)\big)$
for $\theta \ge \lambda^*_+$ and $\nu \in \M$. Consequently,
\begin{equation}\label{eq:F-fdef}
F(\nu):=F(\lambda^*_+,\nu;\Gamma)=\inf_{\theta>\lambda_+^*} \, \{f(\theta,\nu)\} \,, 
\end{equation}
which 
must be in
$[m_-^*,m_+^*]$ 
when $\nu \in \M_1$, 
since then
$f(\theta,\nu)\in[f(\theta,\delta_{\lambda_-^*}),f(\theta,\delta_{\lambda_+^*})]$. 
As $f(\theta,\nu) \to \infty$ when $\theta \to \infty$, 
by monotone convergence, for any $\nu\in \M$ there exists 
some $\theta_\nu \geq \lambda_+^*$ such that $F(\nu)=f(\theta_\nu,\nu)$.
  
\noindent
{\bf I.} \emph{Continuity of $F(\cdot)$}.
The map $\nu \mapsto f(\theta,\nu)$ is continuous on $\M$ 
for each fixed $\theta>\lambda^*_+$, so by \eqref{eq:F-fdef} the infimum 
$F(\nu)$ of these maps is upper semicontinuous (u.s.c.)
on $\M$.  To show that $F(\nu)$ is
lower semi continuous (l.s.c.),
fix a convergent sequence $\nu_n \to \nu$ in $\M$ and let
$\theta_n:=\theta_{\nu_n}$. Passing to a sub-sequence, we may and 
will assume w.l.o.g. that $\theta_n\to \theta^*$ for some $\theta^*$ finite.
If $\theta^*>\lambda_+^*$ then the continuous functions $x\to(\theta_n-x)^{-1}$ 
on the compact $[\lambda_-^*,\lambda_+^*]$
converge uniformly to the continuous function $x\to (\theta^*-x)^{-1}$, from 
which we deduce that as $n \to \infty$, 
$$
F(\nu_n)=f(\theta_n\nu_n) \to f(\theta^*,\nu) \ge F(\nu)
$$
(using \eqref{eq:F-fdef} to get the inequality).
Alternatively, if $\theta_n \to \lambda_+^*$ then $\theta_n \le \lambda_+^*+2 \delta$
for any fixed $\delta>0$ and all $n$ large enough, in which case by monotonicity
of $\theta \mapsto (\theta-x)^{-1}$ and the preceding argument, we have that as
$n \to \infty$,
$$
F(\nu_n) +\delta = f(\theta_n,\nu_n) + \delta \ge f(\lambda_+^*+2\delta,\nu_n)
\to f(\lambda_+^*+2\delta,\nu) \ge F(\nu) \,.
$$ 
Considering $\delta \to 0$ yields the stated l.s.c., hence continuity, of $F(\cdot)$.

\noindent
{\bf II.} \emph{The finiteness of $I^H_q$ and $I^G_q$.}
Setting now $q^{\pm}_t:=t\delta_{\lambda_{\pm}^*}+(1-t) q$, 
both $F(q^+_t) : [0,1] \mapsto [\overline{m},m^*_+]$ and 
$F(q^-_t) : [0,1] \mapsto [m^*_-,\overline{m}]$ are continuous in $t$,
so by the mean-value theorem, for any $m \in (m_-^*,m_+^*)$ there
exists $t=t(m) \in [0,1)$ such that $F(q^{\pm}_t)=m$. Since 
$H(q|q^{\pm}_t)$ are finite, so is $I_q^H(m):=I_q^H(m;\lambda^*_\pm,\Gamma)$.
As for its boundary points, note that $I_q^H(m^*_\pm)=\infty$ unless 
$q=\delta_{\lambda^*_\pm}$, in which case $\overline{m}=m^*_\pm$.
Similarly, $F(y q):\R_+ \mapsto [\frac{1}{2} \lambda^*_+,\infty)$ is 
continuous in $y$, hence for any $m>\frac{1}{2} \lambda^*$ there exists
$y=y(m)>0$ such that $F(y q)=m$. With $H(q|y q)=2 J_1(y)$ finite at
any such $y(m)$, we deduce from the first identity in \eqref{eq:rf-gauss}  
that $I_q^G(m):=I_q^G(m;\lambda^*_\pm,\Gamma)$ is
finite for all $m>\frac{1}{2} \lambda^*_+$.  
 
\noindent
{\bf III.} \emph{Monotonicity and convexity of $I^H_q$ and $I^G_q$.}
Both $I_q^H(m)$ and $I_q^G(m)$ are non-decreasing 
for $m \ge \overline{m}$ in their respective domains. Indeed, 
as seen in step II, there is no need to consider the boundary points.
So, setting either $I_q=I_q^H$ or $I_q=I_q^G$ and fixing $m' \ge m$ 
in the interior of the relevant interval, for any $\epsilon>0$ there 
exists $\nu \in \M$ such that $F(\nu)=m'$ and 
$\frac{1}{2} H(q|\nu) \le I_q(m')+\epsilon$
(with $\nu \in \M_1$ in case $I_q=I_q^H$). 
By the continuity of $F(t \nu + (1-t) q): [0,1] \mapsto [\overline{m},m']$ 
we have that $F(\bar \nu)=m$ for $\bar \nu = s \nu + (1-s) q$ 
and some $s \in [0,1]$. Hence, by the convexity of $\nu\to H(q|\nu)$, 
$$
I_q (m)\leq \frac{1}{2} H(q|\bar\nu)\leq \frac{1}{2} \Big( s H(q|\nu) + (1-s) H(q|q) \Big)
\leq I_q (m')+\epsilon\,.
$$ 
The claimed monotonicity, namely $I_q(m)=\inf_{m' \ge m} \, I_q (m')$,
follows upon considering $\epsilon \downarrow 0$ (and 
by the same reasoning we also get that $I_q (m)$ is 
non-increasing for $m \le \overline{m}$). This monotonicity 
further results with the convexity of $I_q(m)$ for $m \ge \overline{m}$. 
Indeed, for such values of $m$ we have that 
$$
I_q^G (m)=\inf\left\{ \frac{1}{2} H(q|\nu): \nu\in \M, F(\nu) \geq m \right\}\,,
$$
with the analogous formula for $I_q^H(m)$, just requiring then to also have
$\nu \in \M_1$. Now, by the concavity of $F(\cdot)$, if 
$H(q|\nu_i) \le 2 I_q(m_i) + \epsilon$ 
and $F(\nu_i) \ge m_i$, $i=1,2$, then for any $s \in [0,1]$, both 
$F(s \nu_1 + (1-s) \nu_2) \ge s m_1 + (1-s) m_2$ and 
$H(q|s \nu_1 + (1-s) \nu_2) \le 2 (s I_q (m_1) + (1-s) I_q (m_2)) + \epsilon$,
implying that $I_q (s m_1 + (1-s) m_2) \le s I_q (m_1) + (1-s) I_q (m_2)$
(upon taking $\epsilon \downarrow 0$). Finally, since $I_q(m)$ is zero only
at $m=\overline{m}$ and convex at all $m \ge \overline{m}$, it must be 
strictly increasing at any $m \ge \overline{m}$ in its domain. 
 
\noindent  
{\bf IV.} \emph{The continuity of $I^H_q$ and $I^G_q$.}
Clearly the strictly increasing $I_q(m) \to \infty$ when $m \to \infty$.
Thus, the non-negative $I_q(\cdot)$ is a \abbr{GRF} provided it is 
l.s.c. throughout $\R_+$, and to show such l.s.c. it suffices to consider
$m_n \to m$ for which $\alpha := \liminf_{n \to \infty} I_q(m_n)$ is finite.
Recall that $J_1(y) \to \infty$ as $y \to \infty$, so in view of  
\eqref{eq:rf-gauss} we may always restrict our attention to a
compact subset $\M_{[0,y]}=\{\nu \in \M : \nu(\R) \le y\}$ of $\M$
for some $y=y(\alpha) \ge 1$ large enough. Then, by
the continuity of $F(\cdot)$ on the compact $\M_{[0,y]}$ 
we can pass to a sub-sequence $\{n_k\}$ for which there exist 
$\nu_{n_k} \to \nu \in \M_{[0,y]}$ with $F(\nu_{n_k})=m_{n_k}$ 
and $\frac{1}{2} H(q|\nu_{n_k}) \to \alpha$. Since 
$m_{n_k} = F(\nu_{n_k}) \to F(\nu)$  it follows that $F(\nu)=m$ and 
consequently by the l.s.c. of $\nu \mapsto H(q|\nu)$,
$$
\alpha = \frac{1}{2} \lim_{k \to \infty} H(q|\nu_{n_k}) 
\geq \frac{1}{2} H(q|\nu) \ge I_q(m) \,,
$$
as claimed. The continuity of $I_q(\cdot)$ at any $m > \overline{m}$ 
in the interior of its domain, follows from the convexity of $I_q(\cdot)$.
With $m \mapsto I_q(m)$ non-increasing at any $m \le \overline{m}$,
it suffices to fix $\epsilon>0$ and $m \le \overline{m}$ with 
$I_q(m) < \infty$ and show the existence of $m_n \nearrow m$
such that $\liminf_n I_q(m_n) \le I_q(m) + \epsilon$. To this end, 
there exists $\nu \in \M$, $\nu \ne 0$, such that $F(\nu)=m$ and 
$H(q|\nu) \le 2 I_q(m) + \epsilon$ (further having
$\nu \in \M_1$, $\nu \ne \delta_{\lambda^*_-}$ in case $I_q=I_q^H$).
Then, setting $\nu_t = (1-t) \nu + t \delta_{\lambda^*_-} {\bf 1}_{\{I_q=I_q^H\}}$ 
we have that $t \mapsto F(\nu_t)$ is continuous, with 
$H(q|\nu_t) \le H(q|\nu) - \log(1-t)$ and $F(\nu_t)<F(\nu_0)=m$
for all $t \in (0,1]$. Hence, fixing any $t_n \downarrow 0$ results
with $m_n=F(\nu_{t_n}) \nearrow m$, such that 
$$
\liminf_{n \to \infty} I_q(m_n) \le \frac{1}{2} \liminf_{n \to \infty} 
H(q|\nu_{t_n}) \le \frac{1}{2} H(q|\nu) \le I_q(m) + \epsilon
$$
as needed for completing the proof.
\qed

\subsection{A finite dimensional optimization problem}
\label{sec-FD}
For $K>1$ integer, $\bh \in \sqrt{\Gamma} S^{K-1}$, 
$\bl\in \R^K_{\geq}$ and $\bx \in \R^K$, define
$$F_{K,\bh,\bl}(\bx)
=\frac12 \sum_{i=1}^K \lambda_i x_i^2+ \sum_{i=1}^K h_i x_i\,.$$
Our next proposition provides an alternative expression for the optimization problem
\begin{equation}
        \label{eq-Jf}
        F_{K,\bh,\bl}^*=\sup_{\bx\in S^{K-1}} 
        \{F_{K,\bh,\bl}(\bx)\}=\max_{\bx\in S^{K-1}} 
        \{F_{K,\bh,\bl}(\bx)\}\,.
\end{equation}
\begin{proposition}
        \label{prop-JF}
For any $\bh \in \sqrt{\Gamma} S^{K-1}$ and 
$\bl \in \R^K_{\ge}$
let $\nu_{\bh}=\Gamma^{-1} \sum_{i=1}^K h_i^2 \delta_{\lambda_i}$. 
Then,
\begin{equation}
        \label{eq-JF0}
F_{K,\bh,\bl}^*=F(\lambda_1,\nu_{\bh};\Gamma) \,.
        \end{equation}
\end{proposition}

\noindent
Before proving Proposition \ref{prop-JF},
we treat the following easier case.
\begin{lemma}
        \label{lem-Jf}
        Assume $h_1\neq 0$.
        Let $\theta^*$ be the unique solution in 
        $(\lambda_1,\infty)$, of 
        \begin{equation}
                \label{eq-Jftheta}
\sum_{i=1}^K \frac{h_i^2}{(\theta^*-\lambda_i)^2}=\Gamma.
        \end{equation}
Then
\begin{equation}
        \label{eq-JF1}
        F_{K,\bh,\bl}^*=
        \frac12\Big(\theta^* +        
        \sum_{i=1}^K \frac{h_i^2}{\theta^*-\lambda_i}
        \Big)=\frac12
        \inf_{\theta>\lambda_1}
        \Big(\theta +  \sum_{i=1}^K \frac{h_i^2}{\theta-\lambda_i}
        \Big)\,.
\end{equation}
\end{lemma}
\proof[Proof of Lemma \ref{lem-Jf}.]
Note first that $F_{K,\bh,\bl}^*=F_{K,|\bh|,\bl}^*$,
where $|{\bh}|_i=|h_i|$. We thus may assume that $h_i\geq 0$ for all $i$.
By adding a constant to all $\lambda_i$, we may and will also assume that
$\lambda_K>0$.
Finally, with
$\mathcal{B}:=\{{\bf x}: \sum_{i=1}^K  x_i^2\leq 1, 
x_i\geq 0\}$, one 
has from the monotonicity of 
$a\mapsto F_{K,\bh,\bl}(a\bx)$ in $\mathcal{B}$ that
$$\max_{\bx\in \mathcal{B}}
\{F_{K,\bh,\bl}(\bx)\}
=
\sup_{\bx \in S^{K-1}} \{F_{K,\bh,\bl}(\bx)\}.$$
Note also that the maximum of the strictly convex continuous function
$F_{K,\bh,\bl}(\cdot)$ on the convex domain
$\mathcal{B}$
is obtained at a unique 
$\bx^*\in S^{K-1}$ due to the compactness of $S^{K-1}$
and the monotonicity of $a\mapsto F_{K,\bh,\bl}(a\bx)$ in $\mathcal{B}$.

Using the Lagrange multiplier $\frac{\theta}{2}(\sum_{i=1}^K x_i^2 -1)$,
we obtain that $x_i^*(\lambda_i-\theta^*)+ h_i=0$ for all $i$.
This gives $x_i^* =  h_i/(\theta^*-\lambda_i)$ for some $\theta^*$ that
must satisfy \eqref{eq-Jftheta}. Since $x_1^*\geq 0$ is finite and $h_1>0$,
this means in particular that
$\theta^* >  \lambda_1$.
The monotonicity of $\theta\mapsto 
\sum_{i=1}^K h_i^2/(\theta-\lambda_i)^2=:f(\theta)$
on $[\lambda_1,\infty)$ together with $f(\lambda_1)=\infty$, 
        $f(\infty)=0$ yields the uniqueness of such $\theta^*$ satisfying
                \eqref{eq-Jftheta}, as well as the 
                left equality in 
                \eqref{eq-JF1}.
                The second part of \eqref{eq-JF1} then follows by carrying
                out the optimization over $\theta$ 
                in the right hand side and noting that its solution $\bar \theta$ must also satisfy \eqref{eq-Jftheta}, hence coincide with $\theta^*$.
%
                \qed
                \proof[Proof of Proposition \ref{prop-JF}]
                The right side of \eqref{eq-JF1} is precisely $F(\lambda_1,\nu_\bh;\Gamma)$. Hence, in view of Lemma \ref{lem-Jf}, it suffices to consider the 
                case of $h_1=0$, which we handle by approximation. That is, we set $\bh^\epsilon \in \sqrt{\Gamma} S^{K-1}$ so that
                $h_1^\epsilon=\sqrt{\Gamma} \epsilon >0$ and $h_i^\epsilon=\sqrt{1-\epsilon^2} h_i$ for all 
                $i \ge 2$. 
Setting $\phi(\epsilon)=\epsilon + (1-\sqrt{1-\epsilon^2})$, note that
                $$|F^*_{K,\bh^\epsilon,\bl}-F^*_{K,\bh,\bl}|\leq
                 \sqrt{\Gamma} \phi(\epsilon) \to_{\epsilon\to 0} 0\,.$$
Further, 
$\nu_{\bh^\epsilon} \to \nu_{\bh}$ in $\M_1$ hence $F(\lambda_1,\nu_{\bh^\epsilon};\Gamma)
\to F(\lambda_1,\nu_{\bh};\Gamma)$ as
$\epsilon\to 0$ (see part {\bf I} of proof of Proposition \ref{prop:rf-prop}), and
the right side of 
        \eqref{eq-JF1}
        yields \eqref{eq-JF0}. \qed 
        \subsection{An auxiliary LDP for squares of normal variables}
We consider here an auxiliary \abbr{LDP}. Specifically, 
fixing integer $K \ge 1$, partition 
$\{1,\ldots,n\}$ to non-empty, disjoint 
subsets $\I_n(i)$, 
$i=1,\ldots,K$, such that 
$n^{-1} |\I_n(i)| \to_{n\to\infty} \mu_K(i)$,
for $i=1,\ldots,K$, and  
some probability measure $\mu_K$ on $\{1,\ldots,K\}$.
With $\{G_j\}_{j=1}^n$ i.i.d. standard 
normal random variables, define the random 
vectors 
$\bX=\{X_1,\ldots,X_K\}$ 
and $\bbX=\{\bar X_1,\ldots,\bar X_K\}$,
such that $X_i=n^{-1} \sum_{j\in \I_n(i)} G_j^2$ 
and $\bar X_i=X_i/\bX^S$
for $\bX^S=\sum_{i=1}^K X_i$.
Note that the laws of $\bX$ and $\bbX$ depend on $n$.
Finally, let $\S_K=\{\bx\in \R_+^K: \sum_i x_i=1\}$
and associate to each point $\bx\in \S_K$ the 
probability measure $\mu_\bx$ on $\{1,\ldots,K\}$
such that $\mu_\bx(i)=x_i$, $i=1,\ldots,K$. 
\begin{proposition}
        \label{aux-LDP}
        The random vectors $\bbX$ satisfy (as $n\to\infty$) the \abbr{LDP} in
        $\S_K$ with speed $n$ and \abbr{GRF} 
        \begin{equation}
                \label{eq-ratefct}
                J(\bx)=\frac{1}{2} H(\mu_K|\mu_\bx)
        \end{equation}
        where 
        $H(\mu_K|\mu_{\bx})=\sum_{i=1}^K \mu_K(i)\log (\mu_K(i)/\mu_\bx(i))$,
        and we adopt the convention $0\log (0/x)=0$ for all $x\geq 0$.
\end{proposition}
To prove Proposition \ref{aux-LDP}, we first 
establish an elementary result concerning 
large deviations of $\chi$-square variables.
\begin{lemma}
        \label{lem-chisquare}
        Suppose integers $\ell_n \ge 1$ are such that 
        $n^{-1} \ell_n\to \alpha\in [0,1]$. Then,
        $Y_n=n^{-1} \sum_{j=1}^{\ell_n} G_j^2$ 
        satisfies the large deviations on $[0,\infty)$ with \abbr{GRF}
        $J_\alpha (y)=\frac12 (y-\alpha+\alpha\log(\alpha/y))$
        (where again by convention $0 \log (0/x)=0$).
\end{lemma}
\proof A direct computation shows that  
$$\frac1n\log E(e^{\theta n Y_n})=-\frac{\ell_n}{2n} \log (1-2\theta)_+\,.$$
In case $n^{-1} \ell_n \to \alpha>0$ an application of 
the Gartner-Ellis theorem
(see \cite[Theorem 2.3.6]{DZ} for this version),
yields the claim. On the other hand, if $n^{-1} \ell_n \to 0$, fix $y>0$ and 
$\theta_n \uparrow 1/2$ slow enough for 
$n^{-1} \ell_n \log(1-2 \theta_n) \to 0$. Then, 
$$\limsup_{n\to\infty}
\frac1n \log P(Y_n \geq y)\leq -\lim_{n\to\infty}
(\theta_n y + \frac{\ell_n}{2n}\log(1-2\theta_n))
= -\frac{y}{2} = -J_0(y)\,,$$
while since $\ell_n \ge 1$,
$$\liminf_{n \to \infty} \frac1n \log P(Y_n\geq y)\geq \liminf_{n  \to \infty} \frac1n P(G_1^2>ny)=-\frac{y}{2}\,,$$
which completes the proof.
\qed

\proof[Proof of Proposition \ref{aux-LDP}]
Let $f:\R_+^K \mapsto [0,\infty)$ and
$g:\R_+^K\setminus \{0\} \mapsto \S_K$ be defined by 
$f(\bx)=\sum_{i=1}^K x_i$ and 
$g(\bx)=\frac{1}{f(\bx)} \bx$.
By Lemma \ref{lem-chisquare}, and using the independence of its components,
the vector $\bX$ satisfies in $\R_+^K$ the \abbr{LDP} with 
speed $n$ and \abbr{GRF}
$$
\bar J(\bx):=\frac12 \left( f(\bx) - 1 - \log f(\bx) 
+H(\mu_K|\mu_{g(\bx)})\right)\,.
$$
Note that $g(\bX)=\bbX$ and that for any $\delta>0$,
the function $g(\cdot)$ is continuous on 
$f^{-1}((\delta,\infty))$.
Since $\lim_{\delta\to 0} \inf_{\{\bx : f(\bx) \le \delta\}} \, \bar J(\bx)=\infty$, we conclude 
(from the contraction principle,
see \cite[Theorem 4.2.1]{DZ}), that
$\bbX$ satisfies the \abbr{LDP} in $\S_K$ with \abbr{GRF}
$J(\bbx) = \inf_{\{\bx\in \R_+^K: g(\bx)=\bbx\}} \bar J(\bx)$. Clearly, such $J(\bbx)$ is given by 
\eqref{eq-ratefct}, completing the proof.
\qed

\subsection{LDP for quadratic optimization - the diagonal case}\label{sec:diag}

We modify the optimization problem $F^*_{n,\bh,\bl}$ 
so that the \abbr{LDP} of Proposition \ref{aux-LDP} can be applied. 
To this end, for $k=1,2,\ldots$, we let $K=K(k)=2^k$ and form 
refined partitions of the intervals $[\lambda_-^*-\delta_k, \lambda_+^*+\delta_k)$ to disjoint sub-intervals $I_i^{(k)}=[\lambda^-_i,\lambda^+_i)$, 
$i=1,\ldots,K$, such that
$\lambda_i^- < \lambda_+^*$,
$\lambda_i^+ > \lambda_-^*$,
and $q(\{\lambda^\pm_i\})=0$ for $i=1,\ldots,K$, 
while $\Delta_k:=\max_{i=1}^{K} 
(\lambda^+_i-\lambda^-_i)\to 0$ 
as $k\to\infty$ (and with $\Delta_k \ge \delta_k$,
also $\delta_k \to 0$). Let
$\I_n^{(k)}(i) = \{j : \lambda_j(n) \in I_i^{(k)}\}$,
$i=1,\ldots,K$ and 
for any $\bx \in S^{n-1}$ 
set $\bbx \in S^{K-1}$ such that $\bar x_i \ge 0$ and 
$$
\bar x_i^2=\sum_{j \in \I_n^{(k)}(i)} \, x_j^2 \,.
$$
We similarly set $\bbh\in \sqrt{\Gamma}S^{K-1}$ such 
that $\bar h_i\geq 0$ and
$\bar h_i^2=\sum_{j \in \I_n^{(k)} (i)} h_j^2$, 
enforcing $\bar x_i=\bar h_i=0$ in case
the set $\I_n^{(k)}(i)$ is empty. Next, subject to
the latter restriction, define
\begin{eqnarray*}
F^*_{K,\bbh,\bl^\pm} := \sup_{\bbx\in S^{K-1}, \bar{x}_i \geq 0} \Big(\frac12\sum_{i=1}^K \lambda^\pm_i \, \bar x_i^2+\sum_{i=1}^K \bar h_i \, \bar x_i \Big),
\end{eqnarray*}
about the \abbr{LDP} of which we have the following result
(whose proof is deferred to the end of this sub-section).
\begin{proposition}
\label{prop-dQLDP}
Fix $k$ and non-random $\Gamma >0$, taking 
$\bh$ Haar distributed on $\sqrt{\Gamma} S^{n-1}$, independently
of $\bl$.
\newline
(a). The sequence  $\{F^*_{K,\bbh,\bl^-} \}$ 
satisfies the \abbr{LDP}
with speed $n$ and \abbr{GRF}
\begin{equation}\label{eq:rf-haar-k}
I^{H,k}_{q}(m;\Gamma)=
                        \inf \Big\{ \frac12 
                        H(q_k|\nu_k):  
                        m = F(\lambda_1^-,\nu_k;\Gamma) \Big\}\,,
\end{equation}
                        where 
                        $q_k = \sum_{i=1}^K q(I^{(k)}_i) \delta_{\lambda^-_i}$ and
                        $\nu_k = \sum_{i=1}^K 
                        \nu(I^{(k)}_i) \delta_{\lambda^-_i}$
                        for some 
                        $\nu \in M_1(\cup_i I^{(k)}_i)$.     
\newline
(b). For any $m \in \R$, 
\begin{align}\label{eq:conv-of-rf}
I_{q}^H(m;\lambda^*_\pm,\Gamma) &= 
\sup_{\delta>0} \liminf_{k \to \infty} 
\inf_{|m'-m|<\delta} I_{q}^{H,k}(m';\Gamma) \,.  
\end{align}
\end{proposition}

\proof[Proof of Theorem \ref{theo-QLDP}]
$~$ 
\newline
(a). Note that $0\leq F^*_{K,\bbh,\bl^+} - F^*_{K,\bbh,\bl^-} \leq \frac12 \Delta_k$. Further, by Cauchy-Schwarz
\begin{equation}\label{eq:sandw}
F^*_{K,\bbh,\bl^-} \leq F^*_{n,\bh,\bl} \leq F^*_{K,\bbh,\bl^+} \,,
\end{equation}
as soon as $\lambda_n(n) \ge \lambda_-^* - \delta_k$ 
and $\lambda_1(n) \le \lambda_+^* + \delta_k$.
By Assumption \ref{ass-QLDP}, the inequality \eqref{eq:sandw} holds for all 
$n$ large enough, hence
the collection 
$\{F^*_{K,\bbh,\bl^-}\}$ is an 
exponentially good approximation of 
$\{F^*_{n,\bh,\bl}\}$ (see \cite[Definition 4.2.14]{DZ}).
In view of \cite[Theorem 4.2.16, part (a)]{DZ}
(see also \cite[Exercise 4.2.29, part (a)]{DZ}),
part (a) of Theorem \ref{theo-QLDP} is thus a direct consequence of
Proposition \ref{prop-dQLDP}.

\noindent
(b). We represent the 
centered multivariate normal random vector $\bg$  
of covariance matrix $\frac{\Gamma}{n} {\bf I}_n$
as the product of 
Haar distributed $\bh \in \sqrt{\Gamma} S^{n-1}$ 
and the independent $\sqrt{Y_n}$, where $n Y_n$ 
has $\chi$-square law of $n$ degrees of freedom. 
Hence, $F^*_{n,\bg,\bl} = F^*_{n,\sqrt{Y_n}\bh,\bl}$ 
for $Y_n$ of Lemma \ref{lem-chisquare}
(with $\ell_n=n$, so $\alpha=1$),
which is further independent of $\bh$ and $\bl$. 
In particular, the exponentially tight 
$\{Y_n\}$ satisfies the \abbr{LDP}
in $\R_+$ with the \abbr{GRF} $J_1(y)$
of Lemma \ref{lem-chisquare}. 
Moreover, from \eqref{eq-opt} we have that 
$\sqrt{y} \mapsto F^*_{n,\sqrt{y}\bh,\bl}$ is globally Lipschitz
continuous, uniformly in $n$, $\bl$ and $\bh \in \sqrt{\Gamma} S^{n-1}$,
so upon a suitable discretization of the range of $\sqrt{Y_n}$,
we get part (b) of Theorem \ref{theo-QLDP} as an 
immediate consequence of part (a) of this theorem 
(for a similar argument, see \cite[Exercise 4.2.7]{DZ}).
\qed

\proof[Proof of Proposition \ref{prop-dQLDP}]
$~$
\newline
(a). Fixing $k$ and $\Gamma > 0$, we apply Proposition \ref{prop-JF}, 
to find that for each $\bbh\in \sqrt{\Gamma} S^{K-1}$, 
$$
F^*_{K,\bbh,\bl^-}=
F(\lambda_1^-,\nu_{\bbh};\Gamma)
$$
where $\nu_{\bbh}=\Gamma^{-1} \sum_{i=1}^K 
\bar h_i^2 \delta_{\lambda_i^-}$. Next, let
$\cJ_* := \{1,K\} \cup \{1<i<K : q(I_i^{(k)})>0\}$
and note that $\I_n^{(k)}(i)$ is non-empty for 
all $i \in \cJ_*$ and $n \ge n_0(k)$. Indeed, 
with $\lambda^*_-<\lambda_K^+$ and $\lambda_1^- <\lambda^*_+$, we have from (A2) and (A3) that 
both $\I_n^{(k)}(1)$ and $\I_n^{(k)}(K)$ are  
non-empty for all $n \ge n_0(k)$, whereas by
(A1) and our condition that $q(\{\lambda_i^{\pm}\})=0$, the same applies whenever $q(I_i^{(k)})>0$. 
Thus, dividing the positive integers to at 
most $2^{K-2}$ possibilities, we have upon 
passing to the relevant sub-sequence, that 
for some fixed
$\cJ_* \subseteq \cJ \subseteq \{1,\ldots,K\}$
and all $n$, 
$$
n^{-1}|\I_n^{(k)}(i)|=L_n^{\bl}(I^{(k)}_i) > 0
\quad \Longleftrightarrow \quad i \in \cJ \,.
$$ 
Taking $\bh \in \sqrt{\Gamma} S^{n-1}$ according to 
Haar measure, and setting $K'=|\cJ|$, we have that
along such sub-sequence 
$\Gamma^{-1} (\bar h_i^2, i \in \cJ) \in \S_{K'}$
has the law of $\bbX$ of 
Proposition \ref{aux-LDP},
with $\mu_K(i) = \lim_{n \to \infty} L_n^{\bl}(I^{(k)}_i)$
given by $q(I^{(k)}_i)=q_k(\{\lambda^-_i\})$ 
(by Assumption (A1) and having $q(\partial I^{(k)}_i)=0$ for all $i$).  
Now, for any fixed $\Gamma$ and $\{\lambda_i^-, i=1,\ldots,K\}$, the
function $F^*_{K,\bbh,\bl^-}$ of 
$(\bar h_i^2,i \in \cJ)$ is continuous. 
Thus, along such subsequence we get
the \abbr{ldp} in part (a) of Proposition 
\ref{prop-dQLDP} from Proposition \ref{aux-LDP}  
(together with the contraction principle), albeit
having to take in the formula \eqref{eq:rf-haar-k} 
of its \abbr{GRF} only $\nu_k$ supported on
$\cup_{i \in \cJ} I^{(k)}_i$. Further, 
rewriting the proof of Proposition \ref{expl-rf} 
part (a) for $q_k$ and $I_q^{H,k}$
(instead of $q$ and $I^H_q$), we deduce that the 
\abbr{grf} of \eqref{eq:rf-haar-k}
is unchanged by reducing 
the support $\nu_k$, as long as it contains
$\cup_{i \in \cJ_*} I^{(k)}_i$ (see \eqref{eq:rf-gauss-exp}).
This is the case here, regardless of the sub-sequence
we follow, thereby completing the proof of part (a).
\newline
(b). Fixing $\Gamma,m$ 
and turning to the proof of \eqref{eq:conv-of-rf},
note that every $\nu \in M_1([\lambda_-^*,\lambda_+^*])$ of $H(q|\nu)$ finite, 
induces the sequence
$\nu_k = \sum_{i=1}^K \nu(I^{(k)}_i) \delta_{\lambda_i^-}$ 
such that $H(q_k|\nu_k) \uparrow H(q|\nu)$ 
(for example, use $L_1(\nu)$-approximations of 
the relevant bounded continuous test function 
$\phi$ in the variational representation of
\cite[Lemma 6.2.13]{DZ}, by simple functions 
based on the refined partitions $\{I^{(k)}_i\}$).
Further, from \eqref{eq:Mnu-def}
it is easy to see that for any $\xi \ge \lambda_+^*$, 
\begin{equation}\label{eq:M-apx-bds}
F(\xi,\nu;\Gamma) \ge F(\xi,\nu_k;\Gamma) 
\ge 
F(\xi+\Delta_k,\nu;\Gamma)
- \frac12 \Delta_k
\end{equation}
(with the right-inequality holding as soon as
$\xi \ge \lambda_1^-$).
Now, if $I^H_{q}(m;\lambda^*_\pm,\Gamma)<\infty$, then 
for any $\epsilon>0$ there exists
$\nu=\nu^{(\epsilon)} \in M_1([\lambda_-^*,\lambda_+^*])$
such that $\frac12 H(q|\nu) \le I^{H}_q (m;\lambda^*_\pm,\Gamma)+\epsilon$ and $F(\lambda_+^*,\nu;\Gamma) = m$.
Setting $m_k = F(\lambda_1^-,\nu_k;\Gamma)$ for 
$\nu=\nu^{(\epsilon)}$, the latter property
yields, upon considering the left-inequality 
of \eqref{eq:M-apx-bds} 
at $\xi=\lambda_+^* \in [\lambda_1^-,\lambda_1^-+\Delta_k]$ and its 
right-inequality for $\xi=\lambda_1^-$, that
\begin{equation}\label{eq:Mk-cont}
m \ge F(\lambda_+^*,\nu_k;\Gamma) \ge m_k 
\ge m - \frac12 \Delta_k \,.
\end{equation}
With $\Delta_k \to 0$, by \eqref{eq:rf-haar-k} 
and our choice of $\nu=\nu^{(\epsilon)}$,
this implies that for some $m_k \to m$,
$$ 
I^H_q(m;\lambda^*_\pm,\Gamma) + \epsilon \ge 
\frac12 H(q_k|\nu_k) \ge
I^{H,k}_q(m_k;\Gamma) 
\,. 
$$ 
Taking now $\epsilon \downarrow 0$, we conclude that
the l.h.s. of \eqref{eq:conv-of-rf} exceeds its r.h.s. 
\newline
For the converse direction, note that 
$H(q|\nu)=H(q_k|\nu_k)$ for any 
given $k$ and $\nu_k(\cdot)$ considered in 
\eqref{eq:rf-haar-k}, provided 
$\nu \in M_1([\lambda_-^*,\lambda_+^*])$ is given by 
$$
\nu(\cdot):=\sum_{i=1}^K q(\cdot|I_i^{(k)})
\nu_k(\{\lambda^-_i\}) \,.
$$
Further, \eqref{eq:M-apx-bds} holds for
this choice of $\nu(\cdot)$, resulting 
as in the derivation of \eqref{eq:Mk-cont} with
$$
|F(\lambda_+^*,\nu;\Gamma) 
- F(\lambda_1^-,\nu_k;\Gamma)|
\le \frac12  \Delta_k \to 0 \,.
$$
Since this applies \emph{for any} $\nu_k$ which is 
considered in determining $I_{q}^{H,k}(m';\Gamma)$, 
it follows that the r.h.s. of \eqref{eq:conv-of-rf} exceeds 
$$
\sup_{\delta>0} \; \inf_{|m'-m|<2\delta} \{ 
I^{H}_q (m';\lambda^*_\pm,\Gamma) \}
= I^H_q(m;\lambda^*_\pm,\Gamma) 
$$
(due to the lower semi-continuity of 
$I^H_{q}(\cdot;\lambda^*_\pm,\Gamma)$, 
which was proved in Proposition \ref{prop:rf-prop}).
\qed

\begin{remark}\label{rmk:unif}
Denoting by $d_{\textsc{BL}}(\cdot,\cdot)$ the bounded-Lipschitz 
metric compatible with weak convergence in $M_1(\R)$, 
let $B_n((q,\psi^*),\eta)$ denote the collection 
of $\bl \in \R_{\ge}^n$ such that $d_{\textsc{BL}}(L_n^\bl,q)<\eta$
and $|\lambda_1(n)-\psi^*| < \eta$. In conjunction with 
Remark \ref{rem:IG-lminus}, our proof of Theorem \ref{theo-QLDP} 
actually gives for any $m \in \R$ 
the stronger, uniform conclusion in part (b), 
\begin{align}
\lim_{\varepsilon \downarrow 0} &
\liminf_{\eta \downarrow 0} \liminf_{n \to \infty} \frac{1}{n} \log
\inf_{\bl \in B_n((q,\psi^*),\eta)}  P(|F^*_{n,\bg,\bl}-m|<\varepsilon)
= -  I^G_q(m;\psi^*,\Gamma) 
\nonumber \\
= & \lim_{\varepsilon \downarrow 0} \limsup_{\eta \downarrow 0} 
\limsup_{n \to \infty} \frac{1}{n} \log
\sup_{\bl \in B_n((q,\psi^*),\eta)}  P(|F^*_{n,\bg,\bl}-m|<\varepsilon)
 \,,  \qquad \forall \psi^* \ge q_+ \,.
\label{eq:unif-quenched-bds}
\end{align}
The same conclusion applies for the \abbr{LDP} for Haar distributed $\bh$, 
of \abbr{GRF} $I^H_q(m;\lambda^*_\pm,\Gamma)$ which we proved 
in part (a) of Theorem \ref{theo-QLDP}, 
except for replacing in this case $\lambda_1(n)$ by $\lambda_n(n)$ 
whenever $m < \overline{m}$, and considering then $\psi^* \le q_-$.
\end{remark}

\subsection{LDP for matrices: proof of Cor. \ref{cor-QLDP-W} and \ref{cor-ALDP}}

\proof[Proof of Corollary \ref{cor-QLDP-W}]
Fix a sequence of symmetric $\R$-valued 
matrices $\{W_n\} \in \cW_{\lambda^*_\pm,q}$. For each $n$, the matrix $W_n$
of eigenvalue vector $\bl$, 
is of the form $W_n = O_n^T D_n O_n$, for $D_n={\rm diag}(\lambda_1,\ldots,\lambda_n)$ 
and some real, orthogonal matrix $O_n$. Any such
$O_n$ induces the isomorphism  
$\by = O_n \bx$ on $S^{n-1}$, such that 
$\bh=O_n \widetilde{\bh}$ is Haar distributed 
on $\sqrt{\Gamma} S^{n-1}$, independently of $\bl$.
Further, in view of \eqref{eq:Fnw-h-def} and \eqref{eq-opt},
$$
F_n^{W,\widetilde{\bh}}=\sup_{\by\in S^{n-1}} \Big(\frac12 
\langle \bl, \by^2 \rangle + \langle \bh,
\by\rangle\Big) = F^*_{n,\bh,\bl}\,.
$$
Part (a) is thus an immediate consequence of  
part (a) of Theorem \ref{theo-QLDP} 
and the definition of $\cW_{\lambda^*_\pm,q}$. 
Similarly, considering the multivariate normal $\widetilde{\bg}$ 
of covariance $\frac{\Gamma}{n} {\bf I}_n$, results with 
$\bg=O_n \widetilde{\bg}$ having the same law as $\widetilde{\bg}$, independently of $\bl$,
and consequently the \abbr{LDP} of part (b) for $F_n^{W,\widetilde{\bg}}$ follows from
part (b) of Theorem \ref{theo-QLDP} about the \abbr{LDP} of 
$F^*_{n,\bg,\bl}$.
\qed  

\proof[Proof of Corollary \ref{cor-ALDP}] 
We first convert $F^{W,\widetilde{\bg}}_n$ of law $\P_\Gamma^{G,n}$ 
into $F^*_{n,\bg,\bl}$ as in the proof of 
Corollary \ref{cor-QLDP-W}, just now for random $\bl \in \R_{\ge}^n$ 
having the joint eigenvalue density of the \abbr{GOE}.
Recall that the convergence of $L_n^\bl$ to $\sigma$, in $M_1(\R)$, 
occurs with exponential speed $n^2$ (see \cite{BG}). Hence,
fluctuations from this convergence can not affect the 
\abbr{LDP} considered here, which is at exponential speed $n$. 
Specifically, even when proving the \abbr{LDP} upper 
bound, we can assume w.l.o.g. that $d_{\textsc{BL}}(L_n^\bl,\sigma)<\eta$ 
for any $\eta>0$ and all $n \ge n_0(\eta)$. Further,
$\frac{1}{2} \lambda_1 (n) \le F^*_{n,\bg,\bl} \le \frac{1}{2} \lambda_1(n) + \|\bg\|$,
where both $\{\lambda_1(n)\}$ and $\{\|\bg\|\}$ are 
exponentially tight (due to their \abbr{LDP} 
having a \abbr{GRF}, see \cite[Theorem 6.2]{BDG} and
Lemma \ref{lem-chisquare}, respectively). Hence,
the sequence $\{(\lambda_1(n),F^*_{n,\bg,\bl})\}$ is exponentially tight
in $\R^2$, and to establish part (b) of the corollary, it suffices 
to show that for any $m \ge 1$, $\psi^* \ge 2$,
\begin{align}
& \lim_{\varepsilon,\eta \downarrow 0} \,\liminf_{n \to \infty} \;
\frac{1}{n} \log P(\bl \in B_n((\sigma,\psi^*),\eta), 
 \, |F^*_{n,\bg,\bl}-m|<\varepsilon) = 
- [I^G_\sigma (m;\psi^*,\Gamma) + I_e(\psi^*)] 
\nonumber \\
= 
& \lim_{\varepsilon,\eta \downarrow 0} \,\limsup_{n \to \infty} \;
\frac{1}{n} \log P(\bl \in B_n((\sigma,\psi^*),\eta), 
\, |F^*_{n,\bg,\bl}-m|<\varepsilon) 
\label{eq:loc-ldp-ann}
\end{align}
(this is enough due to general considerations, c.f. \cite[Theorems 4.1.11 and 4.2.1]{DZ}).
To this end, fix $m \ge 1$ and $\psi^* \ge 2$. Since the events considered
in \eqref{eq:loc-ldp-ann} are monotone in both $\varepsilon$ and $\eta$, 
we can and will take $\eta \downarrow 0$ before considering $\varepsilon \downarrow 0$. Then, writing 
\begin{align}
 & \frac{1}{n} \log P(\bl \in B_n((\sigma,\psi^*),\eta), 
 \, |F^*_{n,\bg,\bl}-m|<\varepsilon) \nonumber \\
&  = 
\frac{1}{n} \log P \Big(|F^*_{n,\bg,\bl}-m|<\varepsilon \,\Big|\,
\bl \in B_n((\sigma,\psi^*),\eta)\Big) +
\frac{1}{n} \log P\Big( \bl \in B_n((\sigma,\psi^*),\eta) \Big) \,,
\label{eq:decomp-ann}
\end{align}
we have from the uniform bounds of \eqref{eq:unif-quenched-bds}, that 
the term involving
the conditional probability converges to $-I^G_\sigma(m;\psi^*,\Gamma)$
when $n \to \infty$
followed by $\eta \downarrow 0$ and finally $\varepsilon \downarrow 0$.
Further,
due to the much stronger concentration of $L_n^\bl$ under the \abbr{GOE} law, 
for an \abbr{LDP} at exponential speed $n$, 
the events $\{\bl \in B_n((\sigma,\psi^*),\eta)\}$ 
are then equivalent to $\{|\lambda_1(n)-\psi^*|<\eta\}$. Hence, in the limit 
$n \to \infty$ followed by $\eta \downarrow 0$, the right-most term of
\eqref{eq:decomp-ann} converges to $-I_e(\psi^*)$ (by the \abbr{LDP} 
of \cite[Theorem 6.2]{BDG} for 
the top eigenvalue $\{\lambda_1(n)\}$, under the \abbr{GOE} law).
Combining all this, completes the proof of \eqref{eq:loc-ldp-ann} and thereby
of part (b) of the corollary.

Upon replacing $I^G_\sigma(\cdot;\cdot)$ by $I^H_\sigma(\cdot;\cdot)$, the same
argument applies in the Haar setting of $\big\{\P^{H,n}_\Gamma\big\}_{n \ge 1}$
provided $m \ge \overline{m}$. However, to make use of Remark \ref{rmk:unif},
here we must separately consider $m < \overline{m}$, for which the relevant
rare event considered in $B_n((\sigma,\psi^*),\eta)$ is that of having 
$|\lambda_n(n)-\psi^*|<\eta$, for fixed $\psi^* \le -2$ (and all $n$ large enough). 
By the symmetry of the \abbr{GOE} law, the \abbr{LDP} 
for $\{\lambda_n(n)\}$ is up to a sign change of its 
\abbr{GRF}, the same as the \abbr{LDP} for $\{\lambda_1(n)\}$.
Adapting the preceding argument to accomodate for these additional changes, 
takes care of this case as well. The \abbr{grf} $I^{H,A}(m;\Gamma)$ we 
thus obtain for the \abbr{LDP} of $\{\P_\Gamma^{H,n}\}$ matches the 
expression \eqref{eq:ann-rf}, where it is optimal to set $\psi^*_+=2$ 
when $m \le \overline{m}$ and $\psi^*_-=-2$ when $m \ge \overline{m}$. 
\qed

                
\subsection{Rate functions: explicit formulas}

\proof[Proof of Proposition \ref{expl-rf}] 
$~$
\newline
(a). When computing the rate function 
$I^{H}_{q}(\cdot;\lambda^*_\pm,\Gamma)$ 
we consider only $\nu \in M_1([\lambda_-^*,\lambda_+^*])$ 
such that $q \ll \nu$. In particular, 
decomposing such $\nu=\nu_{ac}+\nu_s$ to 
its a.c. and singular parts with respect to $q$, necessarily $\nu_{ac} = \phi q$ for some 
function $\phi$ which is $q$-a.e. positive on 
$[\lambda_-^*,\lambda^*_+]$. Setting 
$t = \nu_s([\lambda_-^*,\lambda^*_+]) = 1-  
\int \phi dq$, elementary algebra 
shows that 
\begin{align}\label{eq:rf-equiv}
\frac12 H(q|\nu) 
& = \int J_1(\phi(x)) q(dx) + \frac{t}{2} \,, \\
F(\lambda_+^*,\nu;\Gamma) & = \frac12 \inf_{\theta > \lambda_+^*} 
\Big[ \theta + \Gamma \int \frac{\phi(x)}{\theta-x} q(dx) + \Gamma \int
\frac{\nu_s(dx)}{\theta-x} \Big]\,,
\label{eq:M-equiv}
\end{align}
with $I_{q}^{H}(m;\lambda_\pm^*,\Gamma)$ 
given by minimizing the r.h.s of 
\eqref{eq:rf-equiv} over non-negative $\phi$ and 
$q$-singular, non-negative measure $\nu_s$ of 
total mass $t=1-\int \phi dq$, subject to the 
given value $m$ of the r.h.s. of \eqref{eq:M-equiv}. 
The r.h.s. of \eqref{eq:rf-equiv} increases 
in $t$, in $(\phi-1)_+$ and in $(1-\phi)_+$, with 
the global minimum (zero) attained at $\phi=1$ 
and $t=0$, for which the expression 
\eqref{eq:M-equiv} equals $\overline{m}$. 
Thus, 
the optimal choice is $\nu_s=t\delta_{\psi^*}$
with $\psi^*=\lambda^*_+$ for $m \ge \overline{m}$ and
$\psi^*=\lambda^*_-$ for $m \le \overline{m}$. 
That is,
\begin{equation}\label{eq:rf-gauss-exp}
I_{q}^{H}(m;\lambda_\pm^*,\Gamma) = 
\inf\left\{ \int J_1(\phi) dq + \frac{t}{2} : 
m = F(\lambda^*_+,\phi q + t \delta_{\psi^*};\Gamma) \right\} \,,
\end{equation}
where we require that $0 \le t = 1 - \int \phi dq$.
Adding to the r.h.s. of \eqref{eq:rf-equiv}
the Lagrange multiplier 
\begin{equation}\label{eq:lagrange}
A [F(\lambda_+^*,\phi q + t \delta_{\psi^*};\Gamma)
-m] + \frac{(B-1)}{2} \big[\int \phi(x) q(dx) + t - 1\big] \,,
\end{equation}
we find that the infimum (over $\phi$), is attained 
for some $\phi^*(x)=(\theta-x)/(B \psi- B x)$ 
(with the equality holding $q$-a.e. 
and $B \psi=B \theta + A \Gamma$). 
Further, per $\phi$ and $t$, the value of 
$F(\lambda_+^*,\phi q+t \delta_{\psi^*};\Gamma)$ 
is attained either at the unique $\theta>\lambda_+^*$ 
for which 
\begin{equation}\label{eq:const0}
D(\theta) := 
\Gamma \int \frac{\phi(x)}{(\theta-x)^2} q(dx) 
+ \frac{\Gamma t}{(\theta-\psi^*)^2}
= 1 \,,
\end{equation}
or at $\theta=\lambda^*_+$, in case 
$D(\lambda^*_+) \le 1$. Now, by our 
assumption that $q_{\pm} = \lambda_{\pm}$, 
the positivity of $\phi^*(\cdot)$ requires 
$\psi \ge \lambda^*_+$, $B>0$, or
$\psi \le \lambda^*_-$, $B<0$ (or 
$B\psi=A\Gamma>0$ when $B=0$), and 
with our Lagrange multiplier we find that
$t=0$ is optimal unless    
$\psi = \psi^*=\lambda^*_+$, $B>0$ 
or $\psi=\psi^*=\lambda^*_-$, $B<0$. 
The constraint $t=1-\int \phi^* dq$ 
amounts to \eqref{eq:const5} and after 
some algebra we deduce that $\phi^*(x)$ 
results with rate function as in \eqref{eq:const4},
where per $m$ (and $B$ satisfying \eqref{eq:const5}),
the values of $\theta,\psi,t$ are determined 
out of \eqref{eq:const2} (the constraint involving 
$m$ in \eqref{eq:rf-gauss-exp}, in case $\phi=\phi^*$), 
and \eqref{eq:const1} 
(which amounts to plugging $\phi=\phi^*$ in \eqref{eq:const0}).
Lastly, as claimed, for $m > \overline{m}$ 
we only consider $\psi^*=\lambda^*_+$, $B>0$ 
and $\psi \in [\psi^*,\theta)$, 
for which $\phi^*(x)$ is increasing 
on $[\lambda^*_-,\lambda^*_+]$, whereas
$m < \overline{m}$ requires $\psi^*=\lambda^*_-$
with either $B>0$, $\psi>\theta$, or $B<0$, 
$\psi \le \psi^*$, in both of which cases  
$\phi^*(x)$ is decreasing on
$[\lambda^*_-,\lambda^*_+]$. 
\newline
(b).
The only difference between $I_q^G(m;\lambda_\pm^*,\Gamma)$ 
and $I_q^H(m;\lambda_\pm^*,\Gamma)$ is that 
any $\nu(\R)>0$ is allowed in the former, so
here $t \ge 0$ and $\int \phi dq \ge 0$ are no
longer constrained to sum to one. Consequently,  
$I^G_q(m;\lambda_\pm^*,\Gamma)$ is also given by 
the r.h.s. of \eqref{eq:rf-gauss-exp}, just 
minimizing now over $\phi \ge 1$, $\psi^* = \lambda^*_+$ and $t \ge 0$ in case $m > \overline{m}$, otherwise 
fixing $t=0$ and minimizing over $\phi \in (0,1]$. 
We proceed as in part (a), except for fixing hereafter $B=1$
in the Lagrange multiplier of \eqref{eq:lagrange}.
Apart from this fixation of $B$, it 
yields the same form of $\phi^*(x)$,  
requiring $t=0$ unless $\psi=\psi^*$ and
having $\theta$ determined by \eqref{eq:const0}.
Also here if $m > \overline{m}$ then we must have
$\psi \in [\lambda_+^*,\theta)$ with 
$t=0$ whenever $\psi>\lambda^*_+$, while 
$\psi \in (\theta,\infty)$ and $t=0$ 
when $m < \overline{m}$. Finally, after some 
algebra we deduce that $\phi^*(x)$ results 
with rate given by \eqref{eq:const3},
for $\theta,\psi,t$ that are determined 
out of \eqref{eq:const2} and \eqref{eq:const1}. 
\qed 

\proof[Proof of Proposition \ref{expl-rf-W}]
$~$
\newline
(a). We are to solve the equations \eqref{eq:const2}--\eqref{eq:const3}
for $B=1$, some $\theta,\psi \ge \psi^*=\lambda^*_+=2$ and 
the semi-circle law $\sigma$. That is, when
$\bG(\xi)= \frac12 [\xi - \sqrt{\xi^2 -4}]$ for $\xi \ge 2$. 
Here $\xi \mapsto 1/\bG(\xi) = \frac12 [\xi + \sqrt{\xi^2-4}]$
is monotone increasing so we can and will change variables to 
$\alpha:=1/\bG(\theta) \ge 1$ and 
$\beta:=1/\bG(\psi) \ge 1$, denoting 
solutions by $(\alpha_q,\beta_q)$. We note that 
$\alpha \ge \beta$ iff $\theta \ge \psi$, which
holds iff $m \ge \overline{m}$, and 
express all quantities appearing in the 
system \eqref{eq:const2}--\eqref{eq:const3}
in terms of $(\alpha,\beta)$. To this end, 
since $\xi=\bG(\xi)+1/\bG(\xi)$ we have that
$\theta=\alpha+1/\alpha$ and $\psi=\beta+1/\beta$.
Further, differentiating we find that 
$d\xi/d\bG(\xi)=1-\bG(\xi)^{-2}$ and hence
\begin{equation}\label{eq:bl-comp}
\bL(\psi)-\bL(\theta)=\int_\theta^\psi \bG(\xi) d\xi =
\int_{1/\alpha}^{1/\beta} 
g(1-g^{-2}) dg = \frac12 (\beta^{-2}-\alpha^{-2}) 
-\log\big(\frac{\alpha}{\beta}\big) \,.
\end{equation}
Combining this with \eqref{eq:const3} yields the formula 
$I^G_\sigma(m;\pm 2,\Gamma)=\I(\alpha,\beta)+t/2$ 
in terms of $\I(\cdot,\cdot)$ of \eqref{eq:T-q-def}.
Turning to determine $(\alpha_q,\beta_q)$ out of 
\eqref{eq:const2} and \eqref{eq:const1}, we have
the following three cases to consider.
\newline
{\bf Case I.} If $\beta_q>1$ then $t=0$
and the unique solution of \eqref{eq:const1} is
$\alpha_q = \max(\beta_q^{-1} (1+\Gamma),1)$. 
Substituting into \eqref{eq:const2}
the option $\alpha_q=1$ results 
with $\beta_q=\frac{\Gamma}{2(m-1)}$. 
However, such a solution can only be relevant if
$$
1 \ge \beta_q^{-1} (1+\Gamma) = 2(m-1)(\Gamma+1)/\Gamma \,,
$$
i.e. for $m \in (1,m_L]$ as in \eqref{eq:theta-is-2}.
\newline
{\bf Case II.} For $\beta_q>1$ and $m>m_L$ we thus must have 
$\beta_q = (1 + \Gamma)/\alpha_q$, which in view of 
\eqref{eq:const2} results with $\alpha_q>1$ such that
$$
2m -\alpha_q^{-1} - \alpha_q = \frac{\Gamma \alpha_q}{1+\Gamma} \,.
$$
This amounts to $\alpha_q>1$ that solve the quadratic equation 
\begin{equation}\label{eq:quad-eq}
m_c^2 \alpha^2 - 2m\alpha + 1 = 0 \,,
\end{equation}
yielding the value of $\alpha_q$ provided in \eqref{eq:int-pt}.
Recall our assumption that the corresponding 
$\beta_q>1$, i.e. that $\alpha_q<1+\Gamma$, which 
for $\alpha_q$ as given in \eqref{eq:int-pt} is 
equivalent to $m \in (m_L,m_U)$. 
\newline
{\bf Case III.} By now we know that for $m \ge m_U$ the only possible solution 
is $\beta=\beta_q=1$ (i.e. $\psi=2$), for which  
\eqref{eq:const2} provides the 
value of $t=\T(\alpha) \ge 0$ as stated in \eqref{eq:T-q-def}.
In this case, upon summing \eqref{eq:const2}   
and \eqref{eq:const1} we deduce that $\alpha=\alpha_q$ must satisfy the 
equality
$$
0 = 2m +\psi-2\theta-\Gamma \bG(\theta)
= 2\big[
m+1 -\alpha - (1+\frac{\Gamma}{2}) \alpha^{-1}\big] \,.
$$ 
The unique $\alpha \ge 1$ that solves this
quadratic equation is given for $m \ge m_U$ by
$\alpha_q$ of \eqref{eq:psi-is-2}. 
\newline
Collecting together {\bf Cases I, II} and {\bf III}, yields the stated formula of 
\eqref{eq:rf-gauss-semi}.

\noindent
(b). Clearly, $I^{G,A}=I^{G}_\sigma$ for all 
$m \le \overline{m}$, since $F^{W,\bh}_n$ is
an increasing function of $\lambda_1(n)$. 
We claim that $I^{G,A}=I^G_\sigma$ also for $m>\overline{m}$,
except when 
setting $B=1$ and $\psi^*>2$ 
in \eqref{eq:const2}-\eqref{eq:const1}, results with $t>0$.  
Indeed, adding the relevant term $I_e(\psi^*)$ 
to the rate function $\int J_1(\phi) dq +t/2$ of 
\eqref{eq:rf-gauss-exp} and using again the
Lagrange multiplier \eqref{eq:lagrange} for $B=1$,
optimality of $\psi^*>2$ requires having 
\begin{equation}\label{eq:psi-star-opt}
I_e'(\psi^*)+\frac{A \Gamma t}{2(\theta-\psi^*)^2} = 0 \,,
\end{equation}
which with $I'_e(\cdot)$ strictly positive, 
implies having $t>0$. Next, recall from the proof of 
part (b) of Proposition \ref{expl-rf} that
$t>0$ requires $\psi^*=\psi=\theta+A\Gamma$.
Since $m>\overline{m}$ we further require 
that $\theta>\psi$ and thus, from \eqref{eq:psi-star-opt} 
deduce that $t=2 (\theta-\psi) I_e'(\psi)$. Plugging 
such value of $t$ into \eqref{eq:const2} and 
\eqref{eq:const1} yields that $\theta>\psi>2$ must be such that 
\begin{align}
2m-\theta &=  \Gamma [\bG(\psi) +  2 I'_e(\psi)] \,,   
\label{eq:const2-ann}
\\
\psi-\theta &= \Gamma [\bG(\theta) -  \bG(\psi) - 2 I'_e(\psi)] \,.
\label{eq:const1-ann}
\end{align}
Next, with $2I'_e(\xi)=\sqrt{\xi^2-4}=1/\bG(\xi)-\bG(\xi)$, 
the identities 
\eqref{eq:const2-ann}-\eqref{eq:const1-ann} are 
in terms of $\alpha=1/\bG(\theta)$ and
$\beta=1/\bG(\psi)$, equivalent to 
\begin{align*}
2m - \alpha^{-1} - \alpha &= \Gamma \beta\,,\\
\beta^{-1} + \beta - \alpha^{-1} - \alpha &= 
\Gamma \alpha^{-1} - \Gamma \beta \,.
\end{align*}
Up to the change $\beta \mapsto \beta^{-1}$, 
these are exactly the equations which determined
$(\alpha_q,\beta_q)$ in {\bf Case II} of part (a).
In conclusion, having $t>0$ requires that we take 
for $\alpha$ the solution $\alpha_a>1$ of the 
quadratic equation \eqref{eq:quad-eq} (which is 
given in \eqref{eq:int-pt}), and then set
$\beta_a^{-1}=\beta_q=(1+\Gamma)/\alpha_a$ for 
the value of $\beta$.
Such solution is only possible if $\beta_a>1$
or equivalently $\alpha_a>1+\Gamma$. As we
have seen before in {\bf Case II} of part (a), 
this amounts to $m>m_U$. 
\newline
Next, similarly to the derivation 
of \eqref{eq:bl-comp}, we find that    
\begin{equation}\label{eq:Ie-comp}
I_e(\psi)=\frac12 \int_{1}^{1/\beta} 
(g^{-1}-g)(1-g^{-2}) dg = 
\frac{1}{4} \big(\beta^2 - \beta^{-2}) - \log \beta \,.
\end{equation}
Further, plugging in \eqref{eq:const3} the optimal 
\begin{equation}\label{eq:opt-t-ann}
t=2 (\theta-\psi) I_e'(\psi) = 
(\alpha + \alpha^{-1} - \beta - \beta^{-1})(\beta-\beta^{-1})\,,
\end{equation}
yields by \eqref{eq:bl-comp} and \eqref{eq:const1-ann}, that for $m > m_U$
\begin{align}
I^G_\sigma(m;\pm\psi^*,\Gamma) &= 
\frac12 [(\theta-\psi) (\bG(\psi) + 2 I_e'(\psi)) 
+ \bL(\psi)-\bL(\theta)] \nonumber \\
&= \frac12 \Big[ 
(\alpha+\alpha^{-1}-\beta-\beta^{-1}) \beta 
+ \frac{1}{2} (\beta^{-2}-\alpha^{-2}) 
-\log\big(\frac{\alpha}{\beta}\big) \Big]\,,
\label{eq:rf-que-ann}
\end{align}
in terms of $\alpha=\alpha_a$ and $\beta=\beta_a$. 
Summing the r.h.s. of \eqref{eq:Ie-comp} and 
\eqref{eq:rf-que-ann}, leads after some algebra 
to the expression $\I(\alpha,\beta^{-1})$. We
have just shown that at $m>m_U$ and
$(\alpha,\beta^{-1})=(\alpha_a,\beta_a^{-1})$ 
given by the r.h.s. of \eqref{eq-addendum},
this is precisely the value of $I^{G,A}(m;\Gamma)$ 
(as stated in \eqref{eq:ann-rf-gauss-exp}).

The function $m \mapsto I^{G,A}(m;\Gamma)$ is clearly smooth everywhere 
except at $m=m_L$. It is further easy to confirm that both $I^{G,A}(m;\Gamma)$
and its first derivative are continuous at $m=m_L$ (where the value of this
function is $\frac{1}{2}(-\log(1-\eta)-\eta-\frac{1}{2}\eta^2)$ and
it derivative equals $-\eta$, for $\eta=\Gamma/(1+\Gamma)$), with 
the second derivative of $I^{G,A}(m;\Gamma)$ being positive everywhere,
thereby verifying its strict convexity.
\qed

\end{document}